\documentclass[10pt, a4paper]{amsart}
\usepackage[dvips,final]{graphics}
\usepackage{array}
\usepackage{arydshln}
\usepackage[makeroom]{cancel}
 \usepackage[all]{xy}
 \usepackage{url}
\usepackage{multirow, blkarray}
\usepackage{booktabs}
\usepackage{textcomp}
\usepackage[final]{epsfig}
\usepackage{color}
\usepackage[T1]{fontenc}      
\usepackage[english,french]{babel}
\usepackage[utf8]{inputenc}
\usepackage{blindtext}

\usepackage{amsfonts,amscd,array, mathdots, epigraph}
\usepackage{amsmath}
\usepackage{amssymb}
\usepackage{amsthm}
\usepackage{mathrsfs}
\usepackage{stmaryrd}
\usepackage{diagbox}
\usepackage{enumitem}
\usepackage{ulem}

\usepackage{ulem}
\usepackage{tikz}
\usepackage{xcolor}
\usepackage{multicol}

\usepackage{hyperref}

\newtheorem{theorem}{Theorem}[section]

\newtheorem{proposition}[theorem]{Proposition}
\newtheorem{corollary}[theorem]{Corollary}

\newtheorem*{examples}{Examples}

\newtheorem{lemma}[theorem]{Lemma}
\newtheorem{definition}[theorem]{Definition}
\newtheorem*{remark}{Remark}

\numberwithin{equation}{section}

\newcommand{\bM}{\mathbb{N^{*}}}

\title{Number of roots of the continuant over a finite local ring}

\author{Flavien Mabilat}

\date{}


\email{flavien.mabilat@univ-reims.fr}

\begin{document}

\maketitle

\selectlanguage{french}
\begin{abstract}

L'objectif de cet article est d'obtenir une formule donnant, pour un entier naturel $n$, le nombre de racines du $n^{i\grave{e}me}$ polynôme continuant lorsque l'on se place sur un anneau local fini. On donnera en particulier des formules permettant de compter les racines du continuant sur les les anneaux locaux $\mathbb{F}_{q}$, $\mathbb{Z}/p^{m}\mathbb{Z}$ et $\frac{\mathbb{F}_{q}[X]}{\langle X^{m} \rangle}$. Pour finir, on donnera, grâce aux méthodes utilisées pour le continuant, une preuve bien plus courte des formules donnant le nombre de $\lambda$-quiddités (qui sont les solutions d'une équation matricielle apparaissant lors de l'étude des frises de Coxeter) sur les anneaux $\mathbb{Z}/p^{m}\mathbb{Z}$.

\end{abstract}

\selectlanguage{english}
\begin{abstract}

The aim of this article is to obtain a formula giving, for a positive integer $n$, the number of roots of the $n^{th}$ continuant polynomial over a finite local ring. In particular, we will give counting formulae for the roots of the continuant over the local rings $\mathbb{F}_{q}$, $\mathbb{Z}/p^{m}\mathbb{Z}$ and $\frac{\mathbb{F}_{q}[X]}{\langle X^{m} \rangle}$. Besides, the methods used for the continuant will allow us to give a new and short proof of the counting formulae for $\lambda$-quiddities (which are the solutions of a matrix equation appearing in the study of Coxeter's friezes) over the rings $\mathbb{Z}/p^{m}\mathbb{Z}$.

\end{abstract}

\selectlanguage{french}

\ \\
\begin{flushright}
\textit{\og C'est encore en méditant l'objet que le sujet a le plus de chance de s'approfondir. \fg} 
\\Gaston Bachelard, \textit{Le Nouvel Esprit scientifique.}
\end{flushright}
\ \\

\selectlanguage{english}

\section{Introduction}
\label{Intro}

Continuants are multivariate polynomials defined recursively over a commutative and unitay ring $A$. More precisely, we set $K_{-1}:=0_{A}$, $K_{0}:=1_{A}$ and for $i \geq 1$ \[K_{i}(X_{1},\ldots,X_{i}):=
\left|
\begin{array}{cccccc}
X_1&1_{A}&&&\\[4pt]
1_{A}&X_{2}&1_{A}&&\\[4pt]
&\ddots&\ddots&\!\!\ddots&\\[4pt]
&&1_{A}&X_{i-1}&\!\!\!\!\!1_{A}\\[4pt]
&&&\!\!\!\!\!1_{A}&\!\!\!\! X_{i}
\end{array}
\right|.\] 

\noindent $K_{i}(X_{1},\ldots,X_{i})$ is the $i^{th}$ continuant and, if there is no ambiguity, we will omit the $i^{th}$ and we will only use the word continuant. To have a concrete expression of these polynomials, one can use Euler's algorithm (see \cite{CO,E}). Besides, for $n \geq 1$, $K_{n}(X_{1},\ldots,X_{n})$ satisfies the two following equalities:
\[K_{n}(X_{1},\ldots,X_{n})=X_{1}K_{n-1}(X_{2},\ldots,X_{n})-K_{n-2}(X_{3},\ldots,X_{n}),\]
\[K_{n}(X_{1},\ldots,X_{n})=X_{n}K_{n-1}(X_{1},\ldots,X_{n-1})-K_{n-2}(X_{1},\ldots,X_{n-2}).\]

These polynomials, which were already known to Euler (see \cite{E}), are involved in the study of many mathematical objects, such as Hirzebruch-Jung continued fractions, discrete Sturm-Liouville equations or Coxeter's friezes (see the introduction of \cite{CO} for more details). In particular, one of the most interesting properties of them is their appearance in the expression of the entries of the following matrix:

\[M_{n}(a_1,\ldots,a_n):=\begin{pmatrix}
    a_{n} & -1_{A} \\[4pt]
     1_{A} & 0_{A}
    \end{pmatrix}
\begin{pmatrix}
    a_{n-1} & -1_{A} \\[4pt]
     1_{A} & 0_{A}
    \end{pmatrix}
    \cdots
    \begin{pmatrix}
    a_{1} & -1_{A} \\[4pt]
     1_{A} & 0_{A}
     \end{pmatrix}.\]
		
\noindent Indeed, for all $n \geq 1$ and for all $(a_1,\ldots,a_n) \in A^{n}$, we have the following equality (see \cite{CO}):
\[M_{n}(a_{1},\ldots,a_{n})=\begin{pmatrix}
    K_{n}(a_{1},\ldots,a_{n}) & -K_{n-1}(a_{2},\ldots,a_{n}) \\
    K_{n-1}(a_{1},\ldots,a_{n-1})  & -K_{n-2}(a_{2},\ldots,a_{n-1}) 
   \end{pmatrix}.\]
		
The matrices $M_{n}(a_{1},\ldots,a_{n})$ play a central role in many topics. In particular, the following equation is very important.
\begin{equation}
\label{p}
\tag{$E_{A}$}
M_{n}(a_1,\ldots,a_n)=\pm Id.
\end{equation}

\noindent Indeed, its solutions are closely linked to the construction of Coxeter's friezes (see for instance \cite{CH}) and to the description of the congruence subgroups of $SL_{2}(\mathbb{Z})$ (see for example \cite{M4,M6}). We will say that a solution $(a_{1},\ldots,a_{n})$ of \eqref{p} is a $\lambda$-quiddity of size $n$ over $A$. There are several ways to study these objects. However, in this article we will only look for counting formulae (see also \cite{BC,CO2,Mo2}).
\\
\\ \indent Let $A$ be a finite commutative and unitary ring. Our main objective here is to find the number of roots of the continuant $K_{n}(X_{1},\ldots,X_{n})$ over $A$. For this, we write, for $n \geq 1$ and $a \in A$: 
\[R_{n}^{a}(A):=\{(a_{1},\ldots,a_{n}) \in A^{n},~K_{n}(a_{1},\ldots,a_{n})=a\},\]
\[r_{n}(A):=\left|R_{n}^{0_{A}}(A)\right|.\]
\noindent Besides, the methods and the results developed for this task will allow us to obtain some counting formulae for $\lambda$-quiddities. Consequently, we will also use the following notations. For $n \geq 1$ and $B \in SL_{2}(A)$, we define the set $\Omega_{n}^{B}(A):=\{(a_{1},\ldots,a_{n}) \in A^{n},~M_{n}(a_{1},\ldots,a_{n})=B\}$ and, if $u$ is an invertible element of $A$ and $B=\begin{pmatrix}
    u & 0_{A} \\
    0_{A}  & u^{-1}
   \end{pmatrix}$, we write:
\[w_{n}^{u}(A):=\left|\Omega_{n}^{B}(A)\right|=\left|\{(a_{1},\ldots,a_{n}) \in A^{n},~M_{n}(a_{1},\ldots,a_{n})=B\}\right|.\]

In all this article, $\mathbb{N}$ is the set of non-negative integers. $(A,+,\times)$ is a commutative and unitary ring different from $\{0_{A}\}$. $0_{A}$ is the identity element of $+$ and $1_{A}$ is the identity element of $\times$. If $l \in \mathbb{N}^{*}$, $l_{A}:=\sum_{i=1}^{l} 1_{A}$. $U(A)$ is the unit group of $A$. If $N \geq 2$ and $a \in \mathbb{Z}$, we set $\overline{a}:=a+N\mathbb{Z}$. Let $N \geq 2$ and $(a,b)\in \mathbb{Z}^{2}$. We write $a \equiv b [N]$ if $N$ divides $a-b$.
\\
\\ \indent Our first aim will be to find an explicit formula for $\sum_{v \in U(A)} w_{n}^{v}(A)$. More precisely, after recalling some useful results on local rings and matrices $M_{n}(a_{1},\ldots,a_{n})$, we will prove, in the section \ref{chap24}, the result below:

\begin{theorem}
\label{prin1}

Let $A$ be a finite local ring and $n \geq 3$. Let $\omega=\sqrt{\left|U(A)\right|^{2}+4 \left|A-U(A)\right|\left|A\right|}$, $\eta=\frac{\left|U(A)\right|+\omega}{2}$ and $\tau=\frac{\left|U(A)\right|-\omega}{2}$. We have the following formula:
\[\sum_{v \in U(A)} w_{n}^{v}(A)=\frac{\left|U(A)\right|\omega(\tau^{n-3}+\eta^{n-3})+(2\left|A-U(A)\right|\left|A\right|+\left|U(A)\right|^{2})(\eta^{n-3}-\tau^{n-3})}{2\omega}.\]

\end{theorem}

In the third section, we will consider the number of roots of the continuant. More precisly, we will prove the two following results.

\begin{theorem}
\label{prin2}

Let $A$ be a finite local ring and $n \geq 1$. 
\[r_{n}(A)=\left|R_{n}^{0_{A}}(A)\right|=\frac{\left|A\right|^{n}-\sum_{u \in U(A)} \left|w_{n+2}^{u}(A)\right|}{\left|A-U(A)\right|}.\]

\noindent Besides, if $n$ is odd, we have: $r_{n}(A)=\frac{\left|A\right|^{n}-\left|U(A)\right|w_{n+2}^{1_{A}}(A)}{\left|A-U(A)\right|}$.

\end{theorem}

If we combine the two results above, we can obtain the values of $r_{n}(A)$ over any finite local rings. For instance, we have:

\begin{corollary}
\label{prin3}

i) Let $n \geq 1$. Let $q$ be a power of a prime number $p$.
 \[r_{n}(\mathbb{F}_{q})=\frac{q^{n}+(-1)^{n+1}}{q+1}.\]

\noindent ii) Let $n \geq 1$, $p$ be a prime number, and $m \geq 2$.
\[r_{n}(\mathbb{Z}/p^{m}\mathbb{Z})=\frac{p^{(m-1)(n-1)}(p^{n}+(-1)^{n+1})}{p+1}.\]

\noindent iii) Let $n \geq 1$, $q$ be a power of a prime number prime number, and $m \geq 1$.
\[r_{n}\left(\frac{\mathbb{F}_{q}[X]}{\langle X^{m} \rangle}\right)=\frac{q^{(m-1)(n-1)}(q^{n}+(-1)^{n+1})}{q+1}.\]

\end{corollary}

In section \ref{chap33}, we will give two different proofs of this corollary. The first one will use theorem \ref{prin1} and the second one will not use it.
\\
\\ \indent Moreover, theorems \ref{prin1} and \ref{prin2} allow us to obtain the number of roots of the continuant over any finite commutative and unitary rings since they are all isomorphic to a direct product of finite local rings (see Theorem \ref{23}).
\\
\\ \indent To conclude, we will give in the last section a new and short proof of the counting formulae for $\lambda$-quiddities over the rings $\mathbb{Z}/p^{m}\mathbb{Z}$ (see \cite{BC} Theorem 1.2). More precisely, we will prove:

\begin{theorem}
\label{prin4}

i) Let $p$ be an odd prime number, $m \geq 2$ and $n=2l \geq 4$.
\begin{itemize}
\item  Suppose that $l$ is even.
\[w_{n+2}^{\overline{1}}(\mathbb{Z}/p^{m}\mathbb{Z})=\frac{p^{(m-1)(n-1)}(p^{l+1}-1)(p^{l}-1)}{p^{2}-1},\]
\[w_{n+2}^{\overline{-1}}(\mathbb{Z}/p^{m}\mathbb{Z})=p^{(m-1)(n-1)}\left(\frac{(p^{l+1}-1)(p^{l}-1)}{p^{2}-1}+p^{l}\right)+\frac{p^{(m-1)(n-1-l)}-1}{p^{n-1-l}-1}p^{ml-1}(p-1).\]
\item Suppose that $l$ is odd.
\[w_{n+2}^{\overline{1}}(\mathbb{Z}/p^{m}\mathbb{Z})=p^{(m-1)(n-1)}\left(\frac{(p^{l+1}-1)(p^{l}-1)}{p^{2}-1}+p^{l}\right)+\frac{p^{(m-1)(n-1-l)}-1}{p^{n-1-l}-1}p^{ml-1}(p-1),\]
\[w_{n+2}^{\overline{-1}}(\mathbb{Z}/p^{m}\mathbb{Z})=\frac{p^{(m-1)(n-1)}(p^{l+1}-1)(p^{l}-1)}{p^{2}-1}.\]
\end{itemize}

\noindent ii) Let $m \geq 3$ and $n=2l \geq 4$.
\begin{itemize}
\item  Suppose that $l$ is even.
\[w_{n+2}^{\overline{1}}(\mathbb{Z}/2^{m}\mathbb{Z})=2^{(m-2)(n-1)} \frac{4^{n}- 2^{n}}{3};~~w_{n+2}^{\overline{-1}}(\mathbb{Z}/2^{m}\mathbb{Z})=2^{(m-2)(n-1)}\frac{4^{n}+2^{n+1}}{3}+\frac{2^{(m-2)(n-1-l)}-1}{2^{n-1-l}-1}2^{ml-1}.\]
\item Suppose that $l$ is odd.
\[w_{n+2}^{\overline{1}}(\mathbb{Z}/2^{m}\mathbb{Z})=2^{(m-2)(n-1)}\frac{4^{n}+ 2^{n+1}}{3}+\frac{2^{(m-2)(n-1-l)}-1}{2^{n-1-l}-1}2^{ml-1};~~w_{n+2}^{\overline{-1}}(\mathbb{Z}/2^{m}\mathbb{Z})=2^{(m-2)(n-1)} \frac{4^{n}-2^{n}}{3}.\]
\end{itemize}

\end{theorem}

\section{Preliminary results}
\label{chap}

\noindent The aim of this section is to give some results which will be useful in the next section. 

\subsection{Local rings}
\label{chap21}

We begin by recalling some basic facts about local rings.

\begin{definition}
\label{21}

Let $A$ be a commutative and unitary ring. $A$ is said to be a local ring if it has a unique maximal ideal.

\end{definition}

\begin{theorem}
\label{21bis}

Let $A$ be a commutative and unitary ring. The following statement are equivalent:
\begin{itemize}
\item $A$ is a local ring;
\item $A-U(A)$ is an ideal;
\item if $a \in A$, then $a$ or $1_{A}-a$ is invertible.
\end{itemize}

\end{theorem}

In particular, if $A$ is a local ring, then for all $a \in A-U(A)$, $1_{A}-a$ and $1_{A}+a$ are invertible (indeed, if $a \in A-U(A)$ then $-a \in A-U(A)$ and $1_{A}+a=1_{A}-(-a) \in U(A)$). 

\begin{proposition}[\cite{BP} proposition 5]
\label{21ter}

The order of a finite local ring is a prime power order.

\end{proposition}

\noindent The next result gives us a list of classical finite local rings.

\begin{proposition}
\label{22}

i) A finite field is a local ring.
\\
\\ii) Let $p$ be a prime number and $m \in \mathbb{N}^{*}$. $\mathbb{Z}/p^{m} \mathbb{Z}$ is a finite local ring.
\\
\\iii) Let $A$ be a finite local ring and $n \in \mathbb{N}^{*}$. $\frac{A[X]}{\langle X^{n} \rangle}$ is a finite local ring.
\\
\\iv) $\frac{(\mathbb{Z}/4\mathbb{Z})[X]}{\langle X^{2}+X+\overline{1} \rangle}$, $\frac{(\mathbb{Z}/4\mathbb{Z})[X]}{\langle X^{2},\overline{2}X \rangle}$ and $\frac{(\mathbb{Z}/2\mathbb{Z})[X,Y]}{\langle X^{2},Y^{2} \rangle}$ are finite local rings.
\\
\\v) A quotient of a local ring is a local ring.
\\
\\vi) A ring isomorphic to a local ring is a local ring.

\end{proposition}

\noindent We can also find commutative and unitary rings which don't satisfy the conditions given in theorem \ref{21bis}.

\begin{examples}
{\rm \begin{itemize}
\item $\mathbb{Z}$ is not a local ring, since, for instance, $1-3=-2\neq \pm 1$.
\item $A=\mathbb{Z}/12\mathbb{Z}$ is not a local ring. Indeed, $\overline{3}$ is a non invertible element of $A$ and $(\overline{1}-\overline{3})=-\overline{2}$ is also a non invertible element of $A$.
\item $A=(\mathbb{Z}/2\mathbb{Z}) \times (\mathbb{Z}/2\mathbb{Z})$ is not a local ring. Indeed, $(\overline{1},\overline{0})$ is a non invertible element of $A$ and $(\overline{1},\overline{1})-(\overline{1},\overline{0})=(\overline{0},\overline{1})$ is also a non invertible element of $A$.
\item $\mathbb{K}$ a commutative field. $A=\left\{\begin{pmatrix}
    x & 0_{\mathbb{K}} \\[4pt]
     0_{\mathbb{K}} & y
    \end{pmatrix}, (x,y) \in \mathbb{K}^{2}\right\}$ is a commutative and unitary ring. However, $\begin{pmatrix}
    1_{\mathbb{K}} & 0_{\mathbb{K}} \\[4pt]
     0_{\mathbb{K}} & 0_{\mathbb{K}}
    \end{pmatrix} \in A-U(A)$ and $\begin{pmatrix}
    1_{\mathbb{K}} & 0_{\mathbb{K}} \\[4pt]
     0_{\mathbb{K}} & 1_{\mathbb{K}}
    \end{pmatrix}-\begin{pmatrix}
    1_{\mathbb{K}} & 0_{\mathbb{K}} \\[4pt]
     0_{\mathbb{K}} & 0_{\mathbb{K}}
    \end{pmatrix}=\begin{pmatrix}
    0_{\mathbb{K}} & 0_{\mathbb{K}} \\[4pt]
     0_{\mathbb{K}} & 1_{\mathbb{K}}
    \end{pmatrix} \in A-U(A)$. So, $A$ is not a local ring.
\end{itemize}
}
\end{examples}

The theorem below shows us that finite local rings are very important in the theory of finite commutative and unitary rings.

\begin{theorem}[\cite{BP} Theorem 1]
\label{23}

Let $A$ be a finite commutative and unitary ring. $A$ is isomorphic to a direct product of finite local rings. If $A \cong A_{1} \times \ldots \times A_{r} \cong B_{1} \times \ldots \times B_{l}$, with $A_{i}, B_{i}$ local, then $r=l$ and there exists a permutation $\sigma$ such that $A_{i} \cong B_{\sigma_{(i)}}$.

\end{theorem}

\subsection{Some basic facts about the matrices $M_{n}(a_{1},\ldots,a_{n})$}
\label{chap22}

The objective of this section is to give some formulae satisfied by the matrices $M_{n}(a_{1},\ldots,a_{n})$ which will be useful in the next section. 

\begin{lemma}
\label{form}

Let $A$ be a commutative and unitary ring and $(a,b,c,u,v) \in A^{5}$.
\\
\\i) $M_{3}(a,1_{A},b)=M_{2}(a-1_{A},b-1_{A})$. 
\\
\\ii) $M_{3}(a,0_{A},b)=-M_{1}(a+b)$.
\\
\\iii) $uv-1_{A}$ invertible. $M_{4}(a,u,v,b)=M_{3}(a+(1_{A}-v)(uv-1_{A})^{-1},uv-1_{A},b+(1_{A}-u)(uv-1_{A})^{-1})$.
\\
\\iv) We suppose that $v$ and $x:=((vb-1_{A})(uv-1_{A})-1_{A})v^{-1}$ are invertible. 
\[M_{5}(a,u,v,b,c)=M_{3}(a-(vb-2_{A})x^{-1},x,c-(uv-2_{A})x^{-1}).\]

\end{lemma}

\begin{proof}

i), ii), and iii) come from \cite{CH} section 4. iv) comes from \cite{BC} (lemma 2.3).

\end{proof}

\begin{lemma}
\label{25}

Let $A$ be a commutative and unitary ring, $(\lambda,u) \in U(A)^{2}$ and $B=\begin{pmatrix}
     u & 0_{A} \\[4pt]
     0_{A} & u^{-1}
    \end{pmatrix}$.
\\
\\i) Let $n=2l \geq 2$. Let $(a_{1},\ldots,a_{n}) \in A^{n}$. If $M_{n}(a_{1},\ldots,a_{n})=B$ then we have the following equality:
\[M_{n}(\lambda a_{1},\lambda^{-1} a_{2},\ldots,\lambda a_{2l-1},\lambda^{-1} a_{2l})=B.\]

\noindent ii) (\cite{BC} lemma 2.5) Let $n=2l+1 \geq 3$. Let $(a_{1},\ldots,a_{n}) \in A^{n}$. If $M_{n}(a_{1},\ldots,a_{n})=B$ then we have the following equality:
\[M_{n}(\lambda a_{1},\lambda^{-1} a_{2},\ldots,\lambda a_{2l-1},\lambda^{-1} a_{2l},\lambda a_{2l+1})=\begin{pmatrix}
     \lambda u & 0_{A} \\[4pt]
     0_{A} & \lambda^{-1} u^{-1}
    \end{pmatrix}.\]
		
\end{lemma}

\begin{proof}

i) comes from a direct adaptation of \cite{CM} lemma 2.16.

\end{proof}

\subsection{Number of $\lambda$-quiddities over a finite local ring}
\label{chap23}

The aim of this subpart is to recall the results which have been already proved. First, we consider the case of finite fields.

\begin{theorem}[Morier-Genoud,~\cite{Mo2} Theorem 1]
\label{26}

Let $q$ be the power of a prime number $p$ and $n>4$.
\\
\\i) $n$ odd. $w_{n}^{-1}(\mathbb{F}_{q})=\frac{q^{n-1}-1}{q^{2}-1}$.
\\
\\ii) If $n$ is even then there exists $m \in \mathbb{N}^{*}$ such that $n=2m$.
\begin{itemize}

\item If $p=2$, $w_{2m}^{-1}(\mathbb{F}_{q})=w_{2m}^{1}(\mathbb{F}_{q})=\frac{(q^{m}-1)(q^{m-1}-1)}{q^{2}-1}+q^{m-1}$.
\item If $p>2$ and $m$ is even, we have: $w_{2m}^{-1}(\mathbb{F}_{q})=\frac{(q^{m}-1)(q^{m-1}-1)}{q^{2}-1}$.
\item If $p>2$ and $m \geq 3$ is odd, we have: $w_{2m}^{-1}(\mathbb{F}_{q})=\frac{(q^{m}-1)(q^{m-1}-1)}{q^{2}-1}+q^{m-1}$.

\end{itemize}

\end{theorem}

\begin{theorem}[\cite{CM} Theorem 1.1]
\label{27}

Let $q$ be the power of a prime number $p>2$ and $n\in\mathbb{N}$, $n>4$.
\\
\\i) If $n$ is odd then we have $w_{n}^{1}(\mathbb{F}_{q})=w_{n}^{-1}(\mathbb{F}_{q})=\frac{q^{n-1}-1}{q^{2}-1}$.
\\
ii) If $n$ is even then there exists $m \in \bM$ such that $n=2m$.
\begin{itemize}
\item If $m$ is even we have: $w_{2m}^{1}(\mathbb{F}_{q})=\frac{(q^{m}-1)(q^{m-1}-1)}{q^{2}-1}+q^{m-1}$.
\item If $m \geq 3$ is odd we have: $w_{2m}^{1}(\mathbb{F}_{q})=\frac{(q^{m}-1)(q^{m-1}-1)}{q^{2}-1}$.
\end{itemize}

\end{theorem}

\begin{remark}
{\rm The previous formulae are still true if $n \in \{2,3,4\}$.
}
\end{remark}

\begin{theorem}[\cite{M5} Theorem 1.4]
\label{29}

Let $m \geq 2$ and $n \geq 2$. We have the following equality:
\[w_{2n+1}^{\overline{1}}(\mathbb{Z}/2^{m}\mathbb{Z})=w_{2n+1}^{\overline{-1}}(\mathbb{Z}/2^{m}\mathbb{Z})=\frac{2^{2mn-2n-2m-1}(2^{2n+3}-8)}{3}.\]

\end{theorem}

\begin{remark}
{\rm The previous formula is still true in the case $m=1$.
}
\end{remark}

\begin{theorem}[\cite{CM} Theorem 1.3]
\label{28}

Let $n=2m \geq 3$ even.

\begin{itemize}
\item If $m$ is even we have $w_{n}^{\overline{1}}(\mathbb{Z}/4\mathbb{Z})=\frac{4^{n-2}+2^{n-1}}{3}$ and $w_{n}^{\overline{-1}}(\mathbb{Z}/4\mathbb{Z})=\frac{4^{n-2}- 2^{n-2}}{3}$.
\item If $m$ is odd we have $w_{n}^{\overline{1}}(\mathbb{Z}/4\mathbb{Z})=\frac{4^{n-2}- 2^{n-2}}{3}$ and $w_{n}^{\overline{-1}}(\mathbb{Z}/4\mathbb{Z})=\frac{4^{n-2}+2^{n-1}}{3}$.
\end{itemize}

\end{theorem}

\begin{theorem}[Böhmler-Cuntz, \cite{BC} Theorem 1.1]
\label{con}

Let $A$ be a finite local ring, $\epsilon \in \{-1,1\}$, and $n>1$ with $n$ odd or $(-1)^{n/2}\neq\epsilon$. Let $K$ be the field $A/(A-U(A))$.
\[w_{n}^{\epsilon 1_{A}}(A)=w_{n}^{\epsilon 1_{K}}(K) \left|A-U(A)\right|^{n-3}.\]

\end{theorem}

If we combine this result with the theorem \ref{23}, we obtain the number of $\lambda$-quiddities of odd size over any finite commutative and unitary ring. Besides, we immediately have the following corollary:

\begin{corollary}
\label{cor}

i) Let $q$ be a power of a prime number, $m \geq 1$, $A=\frac{\mathbb{F}_{q}[X]}{\langle X^{m} \rangle}$ and $n \geq 1$.
\[w_{n+2}^{\overline{1}}(A)=\frac{q^{nm+1}-q^{n(m-1)}}{q^{m-1}(q^{2}-1)}.\]

\noindent ii) (\cite{BC} corollary 3.3) Let $p$ be a prime number, $n \geq 1$ odd and $m \geq 2$.
\[w_{n+2}^{\overline{1}}(\mathbb{Z}/p^{m}\mathbb{Z})=\frac{p^{nm+1}-p^{n(m-1)}}{p^{m-1}(p^{2}-1)}.\]

\end{corollary}

We will also give another proof of i) and two other proofs of ii) which will not use theorem \ref{con} (see the two next sections).

\subsection{Proof of theorem \ref{prin1}}
\label{chap24}

The aim of this section is to prove the counting formula for $\sum_{u \in U(A)} w_{n}^{u}(A)$. The proof given below is an adaptation of the proof of \cite{M5} theorem 1.4 i).

\begin{proof}[Proof of theorem \ref{prin1}]

Let $A$ be a finite local ring, $B \in SL_{2}(A)$ and $n>4$. To simplify our proof, we define the following objects:
\begin{itemize}
\item $\Omega_{n}^{B}(A):=\{(a_{1},\ldots,a_{n}) \in A^{n},~M_{n}(a_{1},\ldots,a_{n})=B\}$;
\item $\Delta_{n}^{B}(A):=\{(a_{1},\ldots,a_{n}) \in \Omega_{n}^{B}(A),~a_{2} \in U(A)\}$;
\item $\Lambda_{n}^{B}(A,a):=\{(a_{1},\ldots,a_{n}) \in \Omega_{n}^{B}(A),~a_{2}=a\}$;
\item $V(A):=\{(u,v) \in A \times U(A),~uv-1_{A} \in A-U(A)\}$;
\item $\begin{array}{ccccc} 
\psi : & V(A) \times A & \longrightarrow & U(A) \\
  & ((u,v),w) & \longmapsto & ((vw-1_{A})(uv-1_{A})-1_{A})v^{-1}  \\
\end{array}$;
\item $T(A,x):=\{((u,v),w) \in V(A) \times A, \psi(u,v,w)=x\}$, $x \in U(A)$.
\end{itemize}

\noindent $\psi$ is well-defined. Indeed, let $(u,v) \in V(A)$ and $w \in A$. The choice of $u$ and $v$ implies $(uv-1_{A}) \notin U(A)$. So, $(vw-1_{A})(uv-1_{A}) \notin U(A)$. Since $A$ is a local ring, $((vw-1_{A})(uv-1_{A})-1_{A}) \in U(A)$. Hence, $((vw-1_{A})(uv-1_{A})-1_{A})v^{-1}$ is invertible. Besides, if $(u,v) \in V(A)$ then $u \in U(A)$ since $A$ is a local ring. We consider the following applications:
\[\begin{array}{ccccc}
\delta : & V(A) & \longrightarrow & U(A) \times (A-U(A)) \\
  & (u,v) & \longmapsto & (u,uv-1_{A}) \\
\end{array}\] and
\[\begin{array}{ccccc}
\tilde{\delta} : & U(A) \times (A-U(A)) & \longrightarrow & V(A) \\
 & (u,x) & \longmapsto & (u,u^{-1}(x+1_{A}))  \\
\end{array}.\]

\noindent $\delta$ and $\tilde{\delta}$ are reciprocal bijections. Hence, $\left|V(A)\right|=\left|U(A)\right|\left|A-U(A)\right|$.
\\
\\ \uwave{First step:} We prove the following recursive formula:
\[\left|\Delta_{n}^{B}(A)\right|=\left|U(A)\right|\left|\Delta_{n-1}^{B}(A)\right|+\left|A-U(A)\right|\left|A\right|\left|\Delta_{n-2}^{B}(A)\right|.\]

\noindent We have the following equalities:
\begin{eqnarray*}
\Delta_{n}^{B}(A) &=& \bigsqcup_{u \in U(A)} \{(a_{1},\ldots,a_{n}) \in \Omega_{n}^{B}(A),~a_{2}=u\} \\
                  &=& \bigsqcup_{u \in U(A)} \{(a_{1},\ldots,a_{n}) \in \Omega_{n}^{B}(A),~a_{2}=u~{\rm and}~ua_{3}-1_{A} \in U(A)\}\\
									& & \bigsqcup_{u \in U(A)} \{(a_{1},\ldots,a_{n}) \in \Omega_{n}^{B}(A),~a_{2}=u~{\rm and}~ua_{3}-1_{A} \in A-U(A)\}\\
						      &=& \bigsqcup_{u \in U(A)} \underbrace{\{(a_{1},\ldots,a_{n}) \in \Omega_{n}^{B}(A),~a_{2}=u~{\rm and}~ua_{3}-1_{A} \in U(A)\}}_{X_{n}^{B}(u)} \\
									& & \bigsqcup_{(u,v) \in V(A)} \underbrace{\{(a_{1},\ldots,a_{n}) \in \Omega_{n}^{B}(A),~a_{2}=u,~a_{3}=v\}}_{Y_{n}^{B}(u,v)}.
\end{eqnarray*}

\noindent Let $u \in U(A)$. First, we consider $X_{n}^{B}(u)$. By lemma \ref{form} iii), we can define the two following applications:

\[\begin{array}{ccccc} 
f_{n,u} : & X_{n}^{B}(u) & \longrightarrow & \Delta_{n-1}^{B}(A) \\
  & (a_{1},u,a_{3},\ldots,a_{n}) & \longmapsto & (a_{1}+(1_{A}-a_{3})(ua_{3}-1_{A})^{-1},ua_{3}-1_{A},a_{4}+(1_{A}-u)(ua_{3}-1_{A})^{-1},a_{5},\ldots,a_{n})  \\
\end{array}\]
\noindent and
\[\begin{array}{ccccc} 
g_{n,u} : & \Delta_{n-1}^{B}(A) & \longrightarrow & X_{n}^{B}(u) \\
 & (a_{1},\ldots,a_{n-1}) & \longmapsto & (a_{1}+(u^{-1}(a_{2}+1_{A})-1_{A})a_{2}^{-1},u,u^{-1}(a_{2}+1_{A}),a_{3}+(u-1_{A})a_{2}^{-1},a_{4},\ldots,a_{n-1})  \\
\end{array}.\]

\noindent $f_{n,u}$ and $g_{n,u}$ are reciprocal bijections. Thus, $\left|X_{n}^{B}(u)\right|=\left|\Delta_{n-1}^{B}(A)\right|$.
\\
\\Now, we set $v \in A$ such that $(u,v) \in V(A)$. We consider $Y_{n}^{B}(u,v)$. We have:
\[Y_{n}^{B}(u,v)=\bigsqcup_{w \in A} \underbrace{\{(a_{1},\ldots,a_{n}) \in \Omega_{n}^{B}(A), a_{2}=u,~a_{3}=v~{\rm and}~a_{4}=w\}}_{Z_{n}^{B}(u,v,w)}.\]
Thus, we will consider the sets $Z_{n}^{B}(u,v,w)$. Let $w \in A$ and $x=\psi(u,v,w)$. By lemma \ref{form} iv), we can define the two following applications:
\[\begin{array}{ccccc} 
h_{n,u,v,w} : & Z_{n}^{B}(u,v,w) & \longrightarrow & \Lambda_{n-2}^{B}(A,x) \\
  & (a_{1},u,v,w,a_{5},\ldots,a_{n}) & \longmapsto & (a_{1}-(vw-2_{A})x^{-1},x,a_{5}-(uv-2_{A})x^{-1},a_{6},\ldots,a_{n})  \\
\end{array}\]
\[{\rm and}~~\begin{array}{ccccc} 
k_{n,u,v,w} : & \Lambda_{n-2}^{B}(A,x) & \longrightarrow & Z_{n}^{B}(u,v,w) \\
 & (a_{1},\ldots,a_{n-2}) & \longmapsto & (a_{1}+(vw-2_{A})x^{-1},u,v,w,a_{3}+(uv-2_{A})x^{-1},a_{4},\ldots,a_{n-2})  \\
\end{array}.\]

\noindent $h_{n,u,v,w}$ and $k_{n,u,v,w}$ are reciprocal bijections. Hence, $\left|Z_{n}^{B}(u,v,w)\right|=\left|\Lambda_{n-2}^{B}(A,x)\right|$.
\\
\\To conclude, we need some information about the sets $T(A,x)$. First, we have: \[V(A) \times A=\bigsqcup_{x \in U(A)} \{((u,v),w) \in V(A) \times A, \psi(u,v,w)=x\}.\]
\noindent Let $x \in U(A)$. We define two applications:
\[\begin{array}{ccccc} 
\alpha_{x} : & T(A,1_{A}) & \longrightarrow & T(A,x) \\
  & ((u,v),w) & \longmapsto & ((ux,vx^{-1}),wx)  \\
\end{array}~~\begin{array}{ccccc}
\beta_{x} : & T(A,x) & \longrightarrow & T(A,1_{A}) \\
 & ((u,v),w) & \longmapsto & ((ux^{-1},vx),wx^{-1})  \\
\end{array}.\]

\noindent $\alpha_{x}$ and $\beta_{x}$ are well-defined and are reciprocal bijections. Hence, $\left|T(A,x)\right|=\left|T(A,1_{A})\right|$. Moreover, 
\[\left|V(A) \times A\right|=\sum_{u \in U(A)} \left|T(A,u)\right|=\left|U(A)\right|\left|T(A,1_{A})\right|.\] 
\noindent So, $\left|T(A,x)\right|=\left|T(A,1_{A})\right|=\left|A-U(A)\right|\left|A\right|$.
\\
\\All these formulae allow us to have the following equalities:
\begin{eqnarray*} 
\left|\Delta_{n}^{B}(A)\right| &=& \sum_{u \in U(A)} \left|X_{n}^{B}(u)\right|+\sum_{(u, v) \in V(A)} \left|Y_{n}^{B}(u,v)\right| \\
                               &=& \sum_{u \in U(A)} \left|\Delta_{n-1}^{B}(A)\right|+\sum_{(u, v) \in V(A)} \left(\sum_{w \in A} \left|Z_{n}^{B}(u,v,w)\right|\right) \\
												       &=& \left|U(A)\right|\left|\Delta_{n-1}^{B}(A)\right|+\sum_{x \in U(A)} \left(\sum_{((u,v),w) \in T(A,x)} \left|Z_{n}^{B}(u,v,w)\right|\right) \\
															 &=& \left|U(A)\right| \left|\Delta_{n-1}^{B}(A)\right|+\sum_{x \in U(A)} \left(\sum_{((u,v),w) \in T(A,x)} \left|\Lambda_{n-2}^{B}(A,x)\right|\right) \\
															 &=& \left|U(A)\right| \left|\Delta_{n-1}^{B}(A)\right|+\sum_{x \in U(A)} \left|T(A,x)\right| \left|\Lambda_{n-2}^{B}(A,x)\right| \\
															 &=& \left|U(A)\right| \left|\Delta_{n-1}^{B}(A)\right|+\left|T(A,1_{A})\right|\sum_{x \in U(A)} \left|\Lambda_{n-2}^{B}(A,x)\right| \\
														   &=& \left|U(A)\right| \left|\Delta_{n-1}^{B}(A)\right|+\left|T(A,1_{A})\right|\left|\Delta_{n-2}^{B}(A)\right| \\
															 &=& \left|U(A)\right| \left|\Delta_{n-1}^{B}(A)\right|+\left|A-U(A)\right|\left|A\right|\left|\Delta_{n-2}^{B}(A)\right|. \\
\end{eqnarray*} 

\noindent \uwave{Second step:} We give a complete formula for $\left|\Delta_{n-1}^{B}(A)\right|$ in the case $B$ diagonal.
\\
\\Let $u \in U(A)$ and $B_{u}:=\begin{pmatrix} 
	 u & 0_{A} \\ 
	 0_{A} & u^{-1}
	\end{pmatrix}$. First, we give a description of $\Delta_{3}^{B_{u}}(A)$ and $\Delta_{4}^{B_{u}}(A)$.
\\
\\i) $M_{3}(a_{1},a_{2},a_{3})=\begin{pmatrix}
    a_{1}a_{2}a_{3}-a_{1}-a_{3} & -a_{2}a_{3}+1_{A} \\[4pt]
     a_{1}a_{2}-1_{A} & -a_{2}
    \end{pmatrix}$. So, $\Delta_{3}^{B_{u}}(A)=\{(-u,-u^{-1},-u)\}$, $\left|\Delta_{3}^{B_{u}}(A)\right|=1$.
\\
\\ii) $M_{4}(a_{1},a_{2},a_{3},a_{4})=\begin{pmatrix}
    a_{1}a_{2}a_{3}a_{4}-a_{3}a_{4}-a_{1}a_{4}-a_{1}a_{2}+1_{A} & -a_{2}a_{3}a_{4}+a_{2}+a_{4} \\[4pt]
     a_{1}a_{2}a_{3}-a_{1}-a_{3} & -a_{2}a_{3}+1_{A}
    \end{pmatrix}$. Hence, we have the following equalities: $\Delta_{4}^{B_{u}}(A)=\{((1-u)v^{-1},v,(1-u^{-1})v^{-1},-uv), v \in U(A)\}$ and $\left|\Delta_{4}^{B_{u}}(A)\right|=\left|U(A)\right|$.
\\
\\Let $\omega:=\sqrt{\left|U(A)\right|^{2}+4\left|A-U(A)\right|\left|A\right|}$, $\eta:=\frac{\left|U(A)\right|+\omega}{2}$ and $\tau:=\frac{\left|U(A)\right|-\omega}{2}$. Let $C:=\begin{pmatrix}
    \left|U(A)\right| & \left|A-U(A)\right|\left|A\right| \\[4pt]
    1    & 0 
   \end{pmatrix}$ and $P:=\begin{pmatrix}
    1 & 1 \\[4pt]
    -\frac{\eta}{\left|A-U(A)\right|\left|A\right|} & -\frac{\tau}{\left|A-U(A)\right|\left|A\right|}
   \end{pmatrix}$. $P$ is invertible and $P^{-1}=\begin{pmatrix}
    \frac{\omega-\left|U(A)\right|}{2\omega} & \frac{-\left|A-U(A)\right|\left|A\right|}{\omega} \\[4pt]
    \frac{\omega+\left|U(A)\right|}{2\omega} & \frac{\left|A-U(A)\right|\left|A\right|}{\omega} 
   \end{pmatrix}$. The characteristic polynomial of $C$ is $\chi_{C}(X)=(X-\eta)(X-\tau)$. We have, $C=P\begin{pmatrix}
    \tau & 0 \\[4pt]
     0 & \eta
   \end{pmatrix}P^{-1}$.
\\		
\\By the first step, we have the following equality:
	\[\begin{pmatrix}
   \left|\Delta_{n}^{B_{u}}(A)\right|  \\[4pt]
    \left|\Delta_{n-1}^{B_{u}}(A)\right|  
   \end{pmatrix}=C\begin{pmatrix}
   \left|\Delta_{n-1}^{B_{u}}(A)\right|  \\[4pt]
    \left|\Delta_{n-2}^{B_{u}}(A)\right|  
   \end{pmatrix}=C^{n-4}\begin{pmatrix}
   \left|\Delta_{4}^{B_{u}}(A)\right|  \\[4pt]
    \left|\Delta_{3}^{B_{u}}(A)\right|  
   \end{pmatrix}=P\begin{pmatrix}
    \tau^{n-4} & 0 \\[4pt]
     0 & \eta^{n-4}
   \end{pmatrix}P^{-1}\begin{pmatrix}
   \left|U(A)\right|  \\[4pt]
    1  
   \end{pmatrix}.\]
	
\noindent Hence, 
\[\left|\Delta_{n}^{B_{u}}(A)\right|=\frac{\left|U(A)\right|\omega(\tau^{n-4}+\eta^{n-4})+(2\left|A-U(A)\right|\left|A\right|+\left|U(A)\right|^{2})(\eta^{n-4}-\tau^{n-4})}{2\omega}.\]

\noindent \uwave{Third step:} We deduce the desired value.
\\
\\Let $v \in U(A)$, $n \geq 2$ and $B_{v}=\begin{pmatrix}
     v & 0_{A} \\[4pt]
     0_{A} & v^{-1}
    \end{pmatrix}$. Lemma \ref{25} ii) and lemma \ref{form} i) allow us to define the two following applications:
 \[\begin{array}{ccccc} 
\vartheta_{n,v} : & \Lambda_{2n+1}^{Id}(A,v) & \longrightarrow & \Omega_{2n}^{B_{v}}(A) \\
  & (a_{1},v,a_{3},\ldots,a_{2n+1}) & \longmapsto & (a_{1}v-1_{A},a_{3}v-1_{A},a_{4}v^{-1},a_{5}v,\ldots,a_{2n+1}v)  \\
\end{array}\]
\[{\rm and}~~\begin{array}{ccccc} 
\theta_{n,v}  : & \Omega_{2n}^{B_{v}}(A) & \longrightarrow & \Lambda_{2n+1}^{Id}(A,v) \\
 & (a_{1},\ldots,a_{2n}) & \longmapsto & ((a_{1}+1_{A})v^{-1},v,(a_{2}+1_{A})v^{-1},a_{3}v,\ldots,a_{2n}v^{-1})  \\
\end{array}.\]

\noindent $\vartheta_{n,v}$ and $\theta_{n,v}$ are reciprocal bijections. Thus, $\left|\Omega_{2n}^{B_{v}}(A)\right|=\left|\Lambda_{2n+1}^{Id}(A,v)\right|$.
\\
\\Hence, $\left|\Delta_{2n+1}^{Id}(A)\right|=\sum_{v \in U(A)} \left|\Lambda_{2n+1}^{Id}(A,v)\right|=\sum_{v \in U(A)} \left|\Omega_{2n}^{B_{v}}(A)\right|=\sum_{v \in U(A)} w_{2n}^{v}(A)$.
\\
\\A similar argument gives $\left|\Omega_{2n+1}^{B_{v}}(A)\right|=\left|\Lambda_{2n+2}^{B_{v}}(A,u)\right|$ for all $u \in U(A)$. So, $\left|U(A)\right|\left|\Omega_{2n+1}^{B_{v}}(A)\right|=\left|\Delta_{2n+2}^{B_{v}}(A)\right|$ and $\sum_{v \in U(A)} w_{2n+1}^{v}(A)=\left|\Delta_{2n+2}^{B_{v}}(A)\right|$ $\left({\rm since}~\left|\Delta_{2n+2}^{B_{v}}(A)\right|~{\rm does~not~depend~on}~v\right)$.
\\
\\If we combine these formulae with the result of the second step, we obtain for all $n \geq 4$: 
\[\sum_{v \in U(A)} w_{n}^{v}(A)=\frac{\left|U(A)\right|\omega(\tau^{n-3}+\eta^{n-3})+(2\left|A-U(A)\right|\left|A\right|+\left|U(A)\right|^{2})(\eta^{n-3}-\tau^{n-3})}{2\omega}.\]

\noindent Besides, this formula is still true if $n=3$.

\end{proof}

\section{Number of roots of the continuant over a finite local ring}
\label{racines}

The aim of this section is to obtain a formula which gives the number of roots of the continuant over a finite local ring. Recall that $R_{n}^{a}(A)=\{(a_{1},\ldots,a_{n}) \in A^{n},~K_{n}(a_{1},\ldots,a_{n})=a\}$.

\subsection{Preliminary results about the roots of the continuant}
\label{chap31}

The aim of this subpart is to collect several results which will be very useful in the proofs of our main results. We begin by the characterization of the roots of $K_{n}(X_{1},\ldots,X_{n})$ for the small values of $n$.

\begin{proposition}
\label{31}

Let $A$ be a commumative and unitary ring.
\\
\\i) $a_{1}$ is a root of $K_{1}(X)$ if and only if $a_{1}=0_{A}$.
\\
\\ii) $(a_{1},a_{2})$ is a root of $K_{2}(X_{1},X_{2})$ if and only if $(a_{1},a_{2}) \in U(A)^{2}$ and $a_{2}=a_{1}^{-1}$. In particular, if $A$ is finite then $K_{2}(X_{1},X_{2})$ has $\left|U(A)\right|$ roots over $A$.

\end{proposition}

\begin{proof}

This is an immediate consequence of the formulae $K_{1}(X)=X$ and $K_{2}(X_{1},X_{2})=X_{1}X_{2}-1_{A}$.

\end{proof}

\begin{lemma}
\label{31bis}

Let $A$ be a commutative and unitary ring. Let $x \in A$, $u \in U(A)$ and $a \in A-U(A)$. 
\\
\\i) Let $n \geq 1$. $R_{n}^{x}(A) \neq \emptyset$.
\\
\\ii) Let $n \geq 2$. There exists $(a_{1},\ldots,a_{n}) \in A^{n}$ such that $K_{n}(a_{1},\ldots,a_{n})=a$ and $K_{n-1}(a_{1},\ldots,a_{n-1})=u$.

\end{lemma}

\begin{proof}

i) By induction on $n$. $K_{1}(x)=x$, $K_{2}(x+1_{A},1_{A})=x$ and so $R_{1}^{x}(A), R_{2}^{x}(A) \neq \emptyset$. Suppose that there exists $n \geq 2$ such that $R_{n}^{x}(A) \neq \emptyset$. There exists $(a_{1},\ldots,a_{n}) \in A^{n}$ such that $K_{n}(a_{1},\ldots,a_{n})=x$. By lemma \ref{form} i), we have $M_{n+1}(a_{1}+1_{A},1_{A},a_{2}+1_{A},\ldots,a_{n})=M_{n}(a_{1},\ldots,a_{n})$. Hence, we have $K_{n+1}(a_{1}+1_{A},1_{A},a_{2}+1_{A},\ldots,a_{n})=K_{n}(a_{1},\ldots,a_{n})=x$ and $R_{n+1}^{x}(A) \neq \emptyset$. By induction, $R_{n}^{x}(A) \neq \emptyset$.
\\
\\ii) By induction on $n$. $K_{2}(u,(a+1_{A})u^{-1})=a$ and $K_{1}(u)=u$. Suppose that there exists $n \geq 2$ such that there exists $(a_{1},\ldots,a_{n}) \in A^{n}$ such that $K_{n}(a_{1},\ldots,a_{n})=a$ and $K_{n-1}(a_{1},\ldots,a_{n-1})=u$. By lemma \ref{form} i), $M_{n+1}(a_{1}+1_{A},1_{A},a_{2}+1_{A},\ldots,a_{n})=M_{n}(a_{1},\ldots,a_{n})$. Hence, $K_{n+1}(a_{1}+1_{A},1_{A},a_{2}+1_{A},\ldots,a_{n})=K_{n}(a_{1},\ldots,a_{n})=a$ and $K_{n}(a_{1}+1_{A},1_{A},a_{2}+1_{A},\ldots,a_{n-1})=K_{n-1}(a_{1},\ldots,a_{n-1})=u$. By induction, the result is true.

\end{proof}

The next result gives us a link between the sets $R_{n}^{a}(A)$ and the elements developed in the previous section.

\begin{lemma}
\label{32}

Let $u \in U(A)$, $B=\begin{pmatrix} 
	u & 0_{A} \\ 
	0_{A} & u^{-1}
	\end{pmatrix}$ and $n \in \mathbb{N}^{*}$. $\left|R_{n}^{u}(A)\right|=\left|\Omega_{n+2}^{-B^{-1}}(A)\right|$.

\end{lemma}

\begin{proof}

We define the following application:
\[\begin{array}{ccccc} 
\Theta_{n,u} : & R_{n}^{u}(A) & \longrightarrow & \Omega_{n+2}^{-B^{-1}}(A) \\
  & (a_{1},\ldots,a_{n}) & \longmapsto & (u^{-1}K_{n-1}(a_{2},\ldots,a_{n}),a_{1},\ldots,a_{n},u^{-1}K_{n-1}(a_{1},\ldots,a_{n-1}))  \\
\end{array}.\]

\noindent Let $(a_{1},\ldots,a_{n}) \in R_{n}^{u}(A)$. We have the following equalities:
\begin{eqnarray*}
k_{1} &:=& K_{n+1}(u^{-1}K_{n-1}(a_{2},\ldots,a_{n}),a_{1},\ldots,a_{n}) \\
      &=& u^{-1}K_{n-1}(a_{2},\ldots,a_{n})K_{n}(a_{1},\ldots,a_{n})-K_{n-1}(a_{2},\ldots,a_{n}) \\
			&=& u^{-1}K_{n-1}(a_{2},\ldots,a_{n})u-K_{n-1}(a_{2},\ldots,a_{n}) \\
			&=& K_{n-1}(a_{2},\ldots,a_{n})-K_{n-1}(a_{2},\ldots,a_{n}) \\
			&=& 0_{A}.
\end{eqnarray*}

\begin{eqnarray*}
k_{2} &:=& K_{n+1}(a_{1},\ldots,a_{n},u^{-1}K_{n-1}(a_{1},\ldots,a_{n-1})) \\
      &=& u^{-1}K_{n-1}(a_{1},\ldots,a_{n-1})K_{n}(a_{1},\ldots,a_{n})-K_{n-1}(a_{1},\ldots,a_{n-1}) \\
			&=& u^{-1}K_{n-1}(a_{1},\ldots,a_{n-1})u-K_{n-1}(a_{1},\ldots,a_{n-1}) \\
			&=& K_{n-1}(a_{1},\ldots,a_{n-1})-K_{n-1}(a_{1},\ldots,a_{n-1}) \\
			&=& 0_{A}.
\end{eqnarray*}

\noindent So, $M:=M_{n+2}(u^{-1}K_{n-1}(a_{2},\ldots,a_{n}),a_{1},\ldots,a_{n},u^{-1}K_{n-1}(a_{1},\ldots,a_{n-1}))=\begin{pmatrix} 
	x & 0_{A} \\ 
	0_{A} & -u
	\end{pmatrix}$, with $x \in A$. Since ${\rm det}(M)=1_{A}$, $x$ is necessarily equal to $-u^{-1}$. Thus, $M=-B^{-1}$ and $\Theta_{n,u}$ is well-defined. Now, we define the following application:
	$\begin{array}{ccccc} 
\tilde{\Theta}_{n,u} : & \Omega_{n+2}^{-B^{-1}}(A) & \longrightarrow & R_{n}^{u}(A) \\
  & (a_{1},\ldots,a_{n+2}) & \longmapsto & (a_{2},\ldots,a_{n+1})  \\
\end{array}$.
\\
\\Let $(a_{1},\ldots,a_{n+2}) \in \Omega_{n+2}^{-B^{-1}}(A)$. 
\[M_{n+2}(a_{1},\ldots,a_{n+2})=\begin{pmatrix}
    K_{n+2}(a_{1},\ldots,a_{n+2}) & -K_{n+1}(a_{2},\ldots,a_{n+2}) \\
    K_{n+1}(a_{1},\ldots,a_{n+1})  & -K_{n}(a_{2},\ldots,a_{n+1}) 
   \end{pmatrix}=\begin{pmatrix}
    -u^{-1} & 0_{A} \\
    0_{A}  & -u 
   \end{pmatrix}.\]
\noindent Hence, $K_{n}(a_{2},\ldots,a_{n+1})=u$ and $\tilde{\Theta}_{n,u}$ is well-defined. 
\\
\\Besides, $\tilde{\Theta}_{n,u}$ and $\Theta_{n,u}$ are reciprocal bijections. Indeed, if $(a_{1},\ldots,a_{n+2}) \in \Omega_{n+2}^{-B^{-1}}(A)$ then we have $0_{A}=a_{1}K_{n}(a_{2},\ldots,a_{n+1})-K_{n-1}(a_{3},\ldots,a_{n+1})$ and $a_{1}=u^{-1}K_{n-1}(a_{3},\ldots,a_{n+1})$. A similar calculation gives $a_{n+2}=u^{-1}K_{n-1}(a_{2},\ldots,a_{n})$.
\\
\\So, $\left|R_{n}^{u}(A)\right|=\left|\Omega_{n+2}^{-B^{-1}}(A)\right|$.

\end{proof}

\begin{lemma}
\label{33}

Let $u \in U(A)$ and $n \in \mathbb{N}^{*}$ odd. $\left|R_{n}^{u}(A)\right|=\left|\Omega_{n+2}^{Id}(A)\right|$.

\end{lemma}

\begin{proof}

Let $n=2l+1$. Let $v=-u^{-1}$ and $B=\begin{pmatrix} 
	v & 0_{A} \\ 
	0_{A} & v^{-1}
	\end{pmatrix}$. We define the two following applications:
\[\begin{array}{ccccc} 
\xi_{n,v} : & \Omega_{n+2}^{Id}(A) & \longrightarrow & \Omega_{n+2}^{B}(A) \\
  & (a_{1},\ldots,a_{2l+3}) & \longmapsto & (v a_{1},v^{-1} a_{2},\ldots,v a_{2l-1},v^{-1} a_{2l+2},v a_{2l+3})  \\
\end{array}\]
\noindent and 
\[\begin{array}{ccccc} 
\zeta_{n,v} : & \Omega_{n+2}^{B}(A) & \longrightarrow & \Omega_{n+2}^{Id}(A) \\
  & (a_{1},\ldots,a_{2l+3}) & \longmapsto & (v^{-1} a_{1},v a_{2},\ldots,v^{-1} a_{2l-1},v a_{2l+2},v^{-1} a_{2l+3})  \\
\end{array}.\]

\noindent Since $n$ is odd, $\xi_{n,v}$ and $\zeta_{n,v}$ are well-defined (lemma \ref{25} ii)). Besides, they are recirocal bijections. Thus, $\left|\Omega_{n+2}^{B}(A)\right|=\left|\Omega_{n+2}^{Id}(A)\right|$. So, by lemma \ref{32}, $\left|R_{n}^{u}(A)\right|=\left|\Omega_{n+2}^{B}(A)\right|=\left|\Omega_{n+2}^{Id}(A)\right|$.

\end{proof}

The following result will be an important tool in the proofs of the counting formulae given in the introduction.

\begin{proposition}
\label{34}

Let $A$ be a finite local ring, $a \in A-U(A)$ and $n \geq 2$. We have the following equality.
\[\left|R_{n}^{a}(A)\right|=\left|R_{n}^{0_{A}}(A)\right|.\]

\end{proposition}

\begin{proof}

Let $n \in \mathbb{N}^{*}$, $n \geq 2$ and $(a_{1},\ldots,a_{n}) \in A^{n}$. Suppose $K_{n}(a_{1},\ldots,a_{n})=x \notin U(A)$. We have
\begin{eqnarray*}
M_{n}(a_{1},\ldots,a_{n}) &=&\begin{pmatrix}
    K_{n}(a_{1},\ldots,a_{n}) & -K_{n-1}(a_{2},\ldots,a_{n}) \\
    K_{n-1}(a_{1},\ldots,a_{n-1})  & -K_{n-2}(a_{2},\ldots,a_{n-1}) 
   \end{pmatrix} \\
	                        &=& \begin{pmatrix}
    x & -K_{n-1}(a_{2},\ldots,a_{n}) \\
    K_{n-1}(a_{1},\ldots,a_{n-1})  & -K_{n-2}(a_{2},\ldots,a_{n-1}) 
   \end{pmatrix}.
\end{eqnarray*}
	
\noindent Hence, $1_{A}={\rm det}(M_{n}(a_{1},\ldots,a_{n}))=-x\times K_{n-2}(a_{2},\ldots,a_{n-1})+K_{n-1}(a_{2},\ldots,a_{n}) K_{n-1}(a_{1},\ldots,a_{n-1})$. So, $K_{n-1}(a_{2},\ldots,a_{n}) K_{n-1}(a_{1},\ldots,a_{n-1})=1_{A}+x\times K_{n-2}(a_{2},\ldots,a_{n-1})$. 
\\
\\Besides, $x\times K_{n-2}(a_{2},\ldots,a_{n-1}) \in A-U(A)$. Indeed, if $x\times K_{n-2}(a_{2},\ldots,a_{n-1}) \in U(A)$ then there exists $u \in U(A)$ such that $1_{A}=(x\times K_{n-2}(a_{2},\ldots,a_{n-1}))u=x(K_{n-2}(a_{2},\ldots,a_{n-1}) \times u)$ and so $x \in U(A)$.
\\
\\Since $A$ is a local ring, $K_{n-1}(a_{2},\ldots,a_{n}) K_{n-1}(a_{1},\ldots,a_{n-1}) \in U(A)$. In particular, $K_{n-1}(a_{1},\ldots,a_{n-1})$ is an invertible element of $A$.
\\
\\Hence, for all $a \in A-U(A)$, we have the following equality: 
\[R_{n}^{a}(A)=\bigsqcup_{u \in U(A)} \underbrace{\{(a_{1},\ldots,a_{n}) \in A^{n}, K_{n}(a_{1},\ldots,a_{n})=a~{\rm and}~K_{n-1}(a_{1},\ldots,a_{n-1})=u\}}_{S_{n}^{a,u}(A)}.\]

\noindent By lemma \ref{31bis} ii), $S_{n}^{a,u}(A) \neq \emptyset$ for all $u \in U(A)$.
\\
\\Let $a \in A-U(A)$ and $u \in U(A)$. We define the two following applications:
\[\begin{array}{ccccc} 
\gamma_{n,a,u} : & S_{n}^{0_{A},u}(A) & \longrightarrow & S_{n}^{a,u}(A) \\
  & (a_{1},\ldots,a_{n}) & \longmapsto & (a_{1},\ldots,a_{n-1},a_{n}+au^{-1})  \\
\end{array},\]
\[\begin{array}{ccccc} 
\kappa_{n,a,u} : & S_{n}^{a,u}(A) & \longrightarrow & S_{n}^{0_{A},u}(A) \\
 & (a_{1},\ldots,a_{n}) & \longmapsto &  (a_{1},\ldots,a_{n-1},a_{n}-au^{-1})  \\
\end{array}.\]

\noindent Let $(a_{1},\ldots,a_{n}) \in S_{n}^{0_{A},u}(A)$. We have the following equalities:
\begin{eqnarray*}
\gamma_{n,a,u}(a_{1},\ldots,a_{n}) &=& K_{n}(a_{1},\ldots,a_{n-1},a_{n}+au^{-1}) \\
                                   &=& (a_{n}+au^{-1})K_{n-1}(a_{1},\ldots,a_{n-1})-K_{n-2}(a_{1},\ldots,a_{n-2}) \\
															     &=& a_{n}K_{n-1}(a_{1},\ldots,a_{n-1})-K_{n-2}(a_{1},\ldots,a_{n-2})+au^{-1}K_{n-1}(a_{1},\ldots,a_{n-1}) \\
															     &=& K_{n}(a_{1},\ldots,a_{n})+au^{-1}u \\
															     &=& 0_{A}+a \\
															     &=& a.
\end{eqnarray*}

\noindent Besides, we obviously have $K_{n-1}(a_{1},\ldots,a_{n-1})=u$. So, $\gamma_{n,a,u}$ is well-defined.
\\
\\Let $(a_{1},\ldots,a_{n}) \in S_{n}^{a,u}(A)$. We have the following equalities:
\begin{eqnarray*}
\kappa_{n,a,u}(a_{1},\ldots,a_{n}) &=& K_{n}(a_{1},\ldots,a_{n-1},a_{n}-au^{-1}) \\
                                   &=& (a_{n}-au^{-1})K_{n-1}(a_{1},\ldots,a_{n-1})-K_{n-2}(a_{1},\ldots,a_{n-2}) \\
															     &=& a_{n}K_{n-1}(a_{1},\ldots,a_{n-1})-K_{n-2}(a_{1},\ldots,a_{n-2})-au^{-1}K_{n-1}(a_{1},\ldots,a_{n-1}) \\
															     &=& K_{n}(a_{1},\ldots,a_{n})-au^{-1}u \\
															     &=& a-a \\
															     &=& 0_{A}.
\end{eqnarray*}

\noindent Besides, we obviously have $K_{n-1}(a_{1},\ldots,a_{n-1})=u$. So, $\kappa_{n,a,u}$ is well-defined.
\\
\\Moreover, $\gamma_{n,a,u}$ and $\kappa_{n,a,u}$ are reciprocal bijections. Hence, $\left|S_{n}^{0_{A},u}(A)\right|=\left|S_{n}^{a,u}(A)\right|$. 
\\
\\So, we have:
\[\left|R_{n}^{a}(A)\right|=\left|\bigsqcup_{u \in U(A)} S_{n}^{a,u}(A)\right|=\sum_{u \in U(A)} \left|S_{n}^{a,u}(A)\right|=\sum_{u \in U(A)} \left|S_{n}^{0_{A},u}(A)\right|=\left|\bigsqcup_{u \in U(A)} S_{n}^{0_{A},u}(A)\right|=\left|R_{n}^{0_{A}}(A)\right|.\]

\end{proof}

\subsection{Proof of theorem \ref{prin2}}
\label{chap32}

The objective of this section is to study the number of roots of the continuant over finite local rings.

\begin{proof}[Proof of theorem \ref{prin2}]

Let $A$ be a finite local ring and $n \geq 2$. We have the following equality:
\begin{eqnarray*}
A^{n} &=& \bigsqcup_{u \in U(A)} \{(a_{1},\ldots,a_{n}) \in A^{n}, K_{n}(a_{1},\ldots,a_{n})=u\} \bigsqcup_{a \in A-U(A)} \{(a_{1},\ldots,a_{n}) \in A^{n}, K_{n}(a_{1},\ldots,a_{n})=a\} \\
      &=& \bigsqcup_{u \in U(A)} R_{n}^{u}(A) \bigsqcup_{a \in A-U(A)} R_{n}^{a}(A).
\end{eqnarray*}

\noindent Besides, by proposition \ref{34}, 
\[\left|\bigsqcup_{a \in A-U(A)} R_{n}^{a}(A)\right|=\sum_{a \in A-U(A)} \left|R_{n}^{a}(A)\right|=\sum_{a \in A-U(A)} \left|R_{n}^{0_{A}}(A)\right|=\left|A-U(A)\right|\left|R_{n}^{0_{A}}(A)\right|.\]

\noindent Hence, $\left|A\right|^{n}=\sum_{u \in U(A)} \left|R_{n}^{u}(A)\right| +\left|A-U(A)\right|\left|R_{n}^{0_{A}}(A)\right|=\sum_{u \in U(A)} w_{n+2}^{-u^{-1}}(A) +\left|A-U(A)\right|r_{n}(A)$ (lemma \ref{32}). Thus, 
\[r_{n}(A)=\frac{\left|A\right|^{n}-\sum_{u \in U(A)} w_{n+2}^{-u^{-1}}(A)}{\left|A-U(A)\right|}=\frac{\left|A\right|^{n}-\sum_{u \in U(A)} w_{n+2}^{u}(A)}{\left|A-U(A)\right|}.\]

\noindent Suppose $n$ odd.
\begin{eqnarray*}
\left|\bigsqcup_{u \in U(A)} R_{n}^{u}(A)\right| &=& \sum_{u \in U(A)} \left|R_{n}^{u}(A)\right| \\
                                                 &=& \sum_{u \in U(A)} \left|\Omega_{n+2}^{Id}(A)\right|~~({\rm lemma}~\ref{33}) \\
																								 &=& \left|U(A)\right| w_{n+2}^{1_{A}}(A).
\end{eqnarray*}

\noindent Hence, $\left|A\right|^{n}=\left|U(A)\right| w_{n+2}^{1_{A}}(A)+\left|A-U(A)\right|\left|R_{n}^{0_{A}}(A)\right|$. Thus, $r_{n}(A)=\left|R_{n}^{0_{A}}(A)\right|=\frac{\left|A\right|^{n}-\left|U(A)\right|w_{n+2}^{1_{A}}(A)}{\left|A-U(A)\right|}$. Besides, this formula is still true if $n=1$.

\end{proof}

\begin{remark}
{\rm There exists infinite sets $R$ (closed under $+$) such that there is a finite number of $\lambda$-quiddities over $R$ of fixed size $n$ and an infinite number of roots of the continuant $K_{n}(X_{1},\ldots,X_{n})$ over $R$. For instance, let $n=5$ and $R=\mathbb{N}^{*}$. There are five $\lambda$-quiddities over $\mathbb{N}^{*}$ (see \cite{CO2}). $K_{n}(X_{1},\ldots,X_{n})$ has infinitely many roots over $\mathbb{N}^{*}$, since for example $(1,1,a,1,1)$ is a root for all $a \in \mathbb{N}^{*}$.
}
\end{remark}

\subsection{Counting formulae over finite fields, over $\mathbb{Z}/p^{m}\mathbb{Z}$ and over $\frac{\mathbb{F}_{q}[X]}{\langle X^{m} \rangle}$}
\label{chap33}

First, we give a proof based on theorems \ref{prin1} and \ref{prin2}.

\begin{proof}[Proof of corollary \ref{prin3}]

i) Let $A=\mathbb{F}_{q}$ with $q$ a prime power number and $n \geq 1$. $A$ is a finite local ring (proposition \ref{22}) and $\left|U(A)\right|=q-1$. $\omega=\sqrt{(q-1)^{2}+4q}=\sqrt{(q+1)^{2}}=q+1$, $\eta=q$ and $\tau=-1$ (notations of theorem \ref{prin1}). Hence, by theorem \ref{prin1},
\begin{eqnarray*}
\sum_{u \in U(A)} w_{n+2,A}^{u} &=& \frac{(q-1)(q+1)((-1)^{n-1}+q^{n-1})+(2q+(q-1)^{2})(q^{n-1}-(-1)^{n-1})}{2(q+1)} \\
                                &=& \frac{(q^{2}-1)(q^{n-1}+(-1)^{n-1})+(q^{2}+1)(q^{n-1}+(-1)^{n})}{2(q+1)} \\
														    &=& \frac{q^{n+1}+(-1)^{n-1}q^{2}-q^{n-1}+(-1)^{n}+q^{n+1}-(-1)^{n-1}q^{2}+q^{n-1}+(-1)^{n}}{2(q+1)} \\
													      &=& \frac{q^{n+1}+(-1)^{n}}{q+1}.
\end{eqnarray*}

\noindent So, by theorem \ref{prin2}, 
\[r_{n}(A)=q^{n}-\frac{q^{n+1}+(-1)^{n}}{q+1}=\frac{q^{n}+(-1)^{n+1}}{q+1}.\]

\noindent ii) Let $A=\mathbb{Z}/p^{m}\mathbb{Z}$ with $p$ a prime number, $m \geq 2$ and $n \geq 1$. $A$ is a finite local ring (proposition \ref{22}) and $\left|U(A)\right|=p^{m-1}(p-1)$. $\omega=\sqrt{p^{2m-2}(p-1)^{2}+4p^{2m-1}}=\sqrt{p^{2m-2}((p-1)^{2}+4p)}=p^{m-1}(p+1)$, $\eta=\frac{p^{m-1}(p-1)+p^{m-1}(p+1)}{2}=p^{m}$ and $\tau=\frac{p^{m-1}(p-1)-p^{m-1}(p+1)}{2}=-p^{m-1}$ (notations of theorem \ref{prin1}). Hence, by theorem \ref{prin1}:
\begin{eqnarray*}
\sum_{u \in U(A)} w_{n+2,A}^{u}  &=& \frac{p^{m-1}(p-1)p^{m-1}(p+1)(\tau^{n-1}+\eta^{n-1})+(2p^{2m-1}+p^{2m-2}(p-1)^{2})(\eta^{n-1}-\tau^{n-1})}{2 p^{m-1}(p+1)} \\
                                 &=& \frac{p^{2m-2}(p^{2}-1)(\tau^{n-1}+\eta^{n-1})+(2p^{2m-1}+p^{2m-2}(p-1)^{2})(\eta^{n-1}-\tau^{n-1})}{2 p^{m-1}(p+1)} \\
                                 &=& \frac{p^{2m-2}(p^{2}-1)(\tau^{n-1}+\eta^{n-1})+(2p^{2m-1}+p^{2m}-2p^{2m-1}+p^{2m-2})(\eta^{n-1}-\tau^{n-1})}{2 p^{m-1}(p+1)} \\
                                 &=& \frac{p^{2m-2}(p^{2}-1)(\tau^{n-1}+\eta^{n-1})+p^{2m-2}(p^{2}+1)(\eta^{n-1}-\tau^{n-1})}{2 p^{m-1}(p+1)} \\
                                 &=& \frac{p^{m+1}\eta^{n-1}-p^{m-1}\tau^{n-1}}{p+1} \\
									               &=& \frac{p^{m+1}p^{m(n-1)}+(-1)^{n}p^{m-1}p^{(n-1)(m-1)}}{p+1} \\
																 &=& \frac{p^{mn+1}+(-1)^{n}p^{(n-1)(m-1)+m-1}}{p+1}.
\end{eqnarray*}

\noindent So, by theorem \ref{prin2}, 
\begin{eqnarray*}
r_{n}(A) &=& \frac{1}{p^{m-1}}\left(p^{mn}-\frac{p^{mn+1}+(-1)^{n}p^{(n-1)(m-1)+m-1}}{p+1}\right) \\
         &=& p^{mn-m+1}-\frac{p^{mn-m+2}+(-1)^{n}p^{(n-1)(m-1)}}{p+1} \\
         &=& \frac{p^{mn-m+2}+p^{nm-m+1}-p^{nm-m+2}+(-1)^{n+1}p^{(n-1)(m-1)}}{p+1} \\
				 &=& \frac{p^{(n-1)(m-1)}(p^{n}+(-1)^{n+1})}{p+1}. 
\end{eqnarray*}																																											
																																										
\noindent iii) Similar calculations give the result.

\end{proof}

Now, we will give another proof of this result which will not use theorem \ref{prin1}. For this, we begin by looking for an explicit formula for $\left|R_{n}^{u}(\mathbb{F}_{q})\right|$ in the cases $u$ invertible.  

\begin{lemma}
\label{35}

Let $q$ be a power of a prime number, $A:=\mathbb{F}_{q}$, $u \in U(A)$ with $u \neq \pm 1_{A}$ and $n=2m \geq 4$. \[\left|R_{2m}^{u}(A)\right|=\frac{(q^{m+1}-1)(q^{m}-1)}{q^{2}-1}.\]

\end{lemma}

\begin{proof} 

Let $A=\mathbb{F}_{q}$, $v=-u^{-1}$, $B=\begin{pmatrix} 
	v & 0_{A} \\ 
	0_{A} & v^{-1}
	\end{pmatrix}$ and $n=2m$. We have the following equality:
\begin{eqnarray*}
\Omega_{n}^{B}(A) &=& \{(a_{1},\ldots,a_{n}) \in \Omega_{n}^{B}(A),~a_{2}=0_{A}\} \sqcup \{(a_{1},\ldots,a_{n}) \in \Omega_{n}^{B}(A),~a_{2} \neq 0_{A}\} \\
                  &=& \bigsqcup_{a \in A} \underbrace{\{(a_{1},\ldots,a_{n}) \in \Omega_{n}^{B}(A),~a_{1}=a~{\rm and}~a_{2}=0_{A}\}}_{\Psi_{a}} \bigsqcup_{b \in A-\{0_{A}\}} \underbrace{\{(a_{1},\ldots,a_{n}) \in \Omega_{n}^{B}(A),~a_{2}=b\}}_{\Theta_{b}} .\\
\end{eqnarray*}
												
\noindent Let $a \in A$. By lemma \ref{form} ii), we can define the two following applications: 
\[\begin{array}{ccccc} 
\varphi_{n,a} : & \Psi_{a} & \longrightarrow & \Omega_{n-2}^{-B}(A) \\
  & (a,0_{A},a_{3},\ldots,a_{2m}) & \longmapsto & (a+a_{3},\ldots,a_{2m})  \\
\end{array}~,~\begin{array}{ccccc} 
\phi_{n,a}  : & \Omega_{n-2}^{-B}(A) & \longrightarrow & \Psi_{a}~~~~~~~~~~~~~~. \\
 & (a_{1},\ldots,a_{2m}) & \longmapsto & (a,0_{A},a_{1}-a,a_{2},\ldots,a_{2m})  \\
\end{array}\]

\noindent $\varphi_{n,a}$ and $\phi_{n,a}$ are reciprocal bijections. Hence, $\left|\Psi_{a}\right|=\left|\Omega_{n-2}^{-B}(A)\right|$.
\\
\\ \noindent Let $b \in A-\{0_{A}\}$. By lemma \ref{25} i), we can define the two following applications: \[\begin{array}{ccccc} 
\vartheta_{n,b} : & \Theta_{b} & \longrightarrow & \Omega_{n-1}^{B}(A) \\
  & (a_{1},b,a_{3},\ldots,a_{2m}) & \longmapsto & (a_{1}b-1_{A},a_{3}b-1_{A},a_{4}b^{-1},a_{5}b,\ldots,a_{2m}b^{-1})  \\
\end{array}\]
\[{\rm and}~~\begin{array}{ccccc} 
\theta_{n,b}  : & \Omega_{n-1}^{B}(A) & \longrightarrow & \Theta_{b} \\
 & (a_{1},\ldots,a_{2m-1}) & \longmapsto & ((a_{1}+1_{A})b^{-1},b,(a_{2}+1_{A})b^{-1},a_{3}b,\ldots,a_{2m-1}b)  \\
\end{array}.\]

\noindent These applications are reciprocal bijections. Thus, $\left|\Theta_{b}\right|=\left|\Omega_{n-1}^{B}(A)\right|$.
\\
\\So, we obtain the following recursive formulae: 
\[w_{2m}^{v}(A)=(q-1)w_{2m-1}^{v}(A)+q w_{2m-2}^{-v}(A)~~ {\rm and}~~w_{2m}^{-v}(A)=(q-1)w_{2m-1}^{-v}(A)+q w_{2m-2}^{v}(A).\]

\noindent Now, we will prove by induction the equality $w_{2l}^{v}(A)=w_{2l}^{-v}(A)$. Since $u \neq \pm 1_{A}$, $v \neq \pm 1_{A}$ and we have $w_{2}^{v}(A)=w_{2}^{-v}(A)=0$. Suppose that there exists $l \geq 1$ such that $w_{2l}^{v}(A)=w_{2l}^{-v}(A)$. By lemma \ref{33}, $w_{2l+1}^{-v}(A)=w_{2l+1}^{v}(A)=w_{2l+1}^{1_{A}}(A)$. So, the two recursive formulae give us $w_{2l+2}^{v}(A)=w_{2l+2}^{-v}(A)$. Hence, by induction, we have, for all $l \geq 1$, $w_{2l}^{v}(A)=w_{2l}^{-v}(A)$. 
\\
\\So, if we use the recursive formula and the theorem \ref{26}, we have the following equalities:
\begin{eqnarray*}
w_{2m}^{v}(A) &=& (q-1)w_{2m-1}^{v}(A)+q w_{2m-2}^{v}(A) \\
						 &=& (q-1)w_{2m-1}^{1_{A}}(A)+q((q-1)w_{2m-3}^{v}(A)+qw_{2m-4}^{v}(A)) \\
						 &=& \ldots~=\sum_{k=0}^{m-2} q^{k}(q-1)w_{2m-(2k+1)}^{1_{A}}(A)+q^{m-1}w_{2}^{v}(A) \\
						 &=& \sum_{k=0}^{m-2} q^{k}(q-1)\frac{q^{2m-2k-2}-1}{q^{2}-1} \\
						 &=& \frac{1}{q+1}\left(q^{2m-2}\sum_{k=0}^{m-2} \frac{1}{q^{k}}-\sum_{k=0}^{m-2} q^{k}\right) \\
						 &=& \frac{1}{q+1}\left(q^{2m-2}\frac{1-q^{1-m}}{1-q^{-1}}-\frac{1-q^{m-1}}{1-q}\right) \\
						 &=& \frac{1}{q+1}\left(q^{2m-1}\frac{1-q^{1-m}}{q-1}-\frac{q^{m-1}-1}{q-1}\right) \\
						 &=& \frac{q^{2m-1}-q^{m}-q^{m-1}+1}{q^{2}-1} \\
						 &=& \frac{(q^{m}-1)(q^{m-1}-1)}{q^{2}-1}. \\
\end{eqnarray*}

\noindent By lemma \ref{32}, $\left|R_{2m}^{u}(A)\right|=w_{2m+2}^{v}(A)=\frac{(q^{m+1}-1)(q^{m}-1)}{q^{2}-1}$.

\end{proof}

\begin{lemma}
\label{36}

Let $n \geq  1$, $m \geq 2$, and $p$ be a prime number. $r_{n}(\mathbb{Z}/p^{m}\mathbb{Z})=p^{n-1}r_{n}(\mathbb{Z}/p^{m-1}\mathbb{Z})$.

\end{lemma}

\begin{proof}

Let $N=p^{m}$. We define the following application:
\[\begin{array}{ccccc} 
\iota : & (\mathbb{Z}/p^{m-1}\mathbb{Z})^{n} & \longrightarrow & \mathcal{P}((\mathbb{Z}/p^{m}\mathbb{Z})^{n}) \\
 & (a_{1}+\frac{N}{p}\mathbb{Z},\ldots,a_{n}+\frac{N}{p}\mathbb{Z}) & \longmapsto & \{(\overline{b_{1}},\ldots,\overline{b_{n}}),\exists k_{i} \in [\![0;p-1]\!]~\overline{b_{i}}=\overline{a_{i}+\frac{k_{i}N}{p}}\}  \\
\end{array}.\]

\noindent Our first aim is to prove $\bigsqcup_{z \in R_{n}^{0}(\mathbb{Z}/p^{m-1}\mathbb{Z})} \iota(z)= \bigsqcup_{k=0}^{p-1} R_{n}^{\overline{\frac{kN}{p}}}(\mathbb{Z}/N \mathbb{Z})$.
\\
\\Let $(a_{1},\ldots,a_{n},c_{1},\ldots,c_{n}) \in \mathbb{Z}^{2n}$ such that $z:=(a_{1}+\frac{N}{p}\mathbb{Z},\ldots,a_{n}+\frac{N}{p}\mathbb{Z}) \neq z':=(c_{1}+\frac{N}{p}\mathbb{Z},\ldots,c_{n}+\frac{N}{p}\mathbb{Z})$. Suppose $\iota(z) \cap \iota(z') \neq \emptyset$. There exists $1 \leq j \leq n$ such that $a_{j}+\frac{N}{p}\mathbb{Z} \neq c_{j}+\frac{N}{p}\mathbb{Z}$ and there exists $(\overline{b_{1}},\ldots,\overline{b_{n}}) \in \iota(z) \cap \iota(z')$. So, there exists $(k,l) \in [\![0;p-1]\!]^{2}$ such that $\overline{b_{j}}=\overline{a_{j}+\frac{kN}{p}}$ and $\overline{b_{j}}=\overline{c_{j}+\frac{lN}{p}}$. Thus, $a_{j}+\frac{kN}{p}+N\mathbb{Z}=c_{j}+\frac{lN}{p}+N\mathbb{Z}$ and $a_{j}+\frac{N}{p}\mathbb{Z}=c_{j}+\frac{N}{p}\mathbb{Z}$. Thus, $\iota(z) \cap \iota(z')=\emptyset$. 
\\
\\Let $(a_{1},\ldots,a_{n}) \in \mathbb{Z}^{n}$ and $(\epsilon_{1},\ldots,\epsilon_{n}) \in [\![0;p-1]\!]^{n}$. $K_{n}(a_{1}+\epsilon_{1} \frac{N}{p},\ldots,a_{n}+\epsilon_{n}\frac{N}{p}) \equiv K_{n}(a_{1},\ldots,a_{n}) \left[\frac{N}{p}\right]$.
\\
\\Suppose $(a_{1}+\frac{N}{p}\mathbb{Z},\ldots,a_{n}+\frac{N}{p}\mathbb{Z}) \in R_{n}^{0}(\mathbb{Z}/p^{m-1}\mathbb{Z})$ then $K_{n}(a_{1}+\epsilon_{1} \frac{N}{p},\ldots,a_{n}+\epsilon_{n}\frac{N}{p}) \equiv 0 \left[\frac{N}{p}\right]$. Thus, there exists $k \in [\![0;p-1]\!]$ such that $K_{n}(a_{1}+\epsilon_{1} \frac{N}{p},\ldots,a_{n}+\epsilon_{n}\frac{N}{p}) \equiv \frac{kN}{p} [N]$. Hence, $(\overline{a_{1}+\epsilon_{1} \frac{N}{p}},\ldots,\overline{a_{n}+\epsilon_{n}\frac{N}{p}}) \in \bigsqcup_{k=0}^{p-1} R_{n}^{\overline{\frac{kN}{p}}}(\mathbb{Z}/N \mathbb{Z})$, that is to say $\iota\left(\left(a_{1}+\frac{N}{p}\mathbb{Z},\ldots,a_{n}+\frac{N}{p}\mathbb{Z}\right)\right) \subset \bigsqcup_{k=0}^{p-1} R_{n}^{\overline{\frac{kN}{p}}}(\mathbb{Z}/N \mathbb{Z})$. 
\\
\\Suppose $(\overline{a_{1}},\ldots,\overline{a_{n}}) \in \bigsqcup_{k=0}^{p-1} R_{n}^{\overline{\frac{kN}{p}}}(\mathbb{Z}/N \mathbb{Z})$. There exists $k \in [\![0;p-1]\!]$ such that $K_{n}(a_{1},\ldots,a_{n}) \equiv \frac{kN}{p} [N]$. So, $K_{n}(a_{1},\ldots,a_{n}) \equiv 0 \left[\frac{N}{p}\right]$. Thus, $(a_{1}+\frac{N}{p}\mathbb{Z},\ldots,a_{n}+\frac{N}{p}\mathbb{Z}) \in R_{n}^{0}(\mathbb{Z}/p^{m-1}\mathbb{Z})$ and $(\overline{a_{1}},\ldots,\overline{a_{n}}) \in \iota((a_{1}+\frac{N}{p}\mathbb{Z},\ldots,a_{n}+\frac{N}{p}\mathbb{Z}))$.
\\
\\So, $\bigsqcup_{z \in R_{n}^{0}(\mathbb{Z}/p^{m-1}\mathbb{Z})} \iota(z)=\bigsqcup_{k=0}^{p-1} R_{n}^{\overline{\frac{kN}{p}}}(\mathbb{Z}/N \mathbb{Z})$.
\\
\\Besides, $\mathbb{Z}/N \mathbb{Z}$ is a local ring and for all $k \in [\![0;p-1]\!]$ $\overline{\frac{kN}{p}}$ is non invertible (since $m \geq 2$). So, by proposition \ref{34}, we have for all $k \in [\![0;p-1]\!]$:
\[\left|R_{n}^{\overline{0}}(\mathbb{Z}/N \mathbb{Z})\right|=\left|R_{n}^{\overline{\frac{kN}{p}}}(\mathbb{Z}/N \mathbb{Z})\right|.\] 

\noindent If we collect all the preceding elements, we have:
\begin{eqnarray*}
p~r_{n}(\mathbb{Z}/p^{m}\mathbb{Z})=p\left|R_{n}^{\overline{0}}(\mathbb{Z}/N \mathbb{Z})\right| &=& \sum_{k=0}^{p-1} \left|R_{n}^{\overline{\frac{kN}{p}}}(\mathbb{Z}/N \mathbb{Z})\right| \\
                                                                                                &=& \left|\bigsqcup_{k=0}^{p-1} R_{n}^{\overline{\frac{kN}{p}}}(\mathbb{Z}/N \mathbb{Z})\right| \\
																		                                                            &=& \left|\bigsqcup_{z \in R_{n}^{0}(\mathbb{Z}/p^{m-1}\mathbb{Z})} \iota(z)\right| \\
                                                                                                &=& \sum_{z \in R_{n}^{0}(\mathbb{Z}/p^{m-1}\mathbb{Z})} \left|\iota(z)\right| \\
                                                                                                &=& \sum_{z \in R_{n}^{0}(\mathbb{Z}/p^{m-1}\mathbb{Z})} p^{n} \\
																		                                                            &=& p^{n}r_{n}(\mathbb{Z}/p^{m-1}\mathbb{Z}).
\end{eqnarray*}

\noindent Hence, $r_{n}(\mathbb{Z}/p^{m}\mathbb{Z})=p^{n-1}r_{n}(\mathbb{Z}/p^{m-1}\mathbb{Z})$.

\end{proof}

\begin{lemma}
\label{37}

Let $n \geq  1$, $m \geq 2$, and $q$ be a power of a prime number. $r_{n}\left(\frac{\mathbb{F}_{q}[X]}{\langle X^{m} \rangle}\right)=q^{n-1}r_{n}\left(\frac{\mathbb{F}_{q}[X]}{\langle X^{m-1} \rangle}\right)$.

\end{lemma}

\begin{proof}

The arguments are overall the same as those used to prove the previous proposition. We define the following application:
\[\begin{array}{ccccc} 
\tilde{\iota} : & \left(\frac{\mathbb{F}_{q}[X]}{\langle X^{m-1} \rangle}\right)^{n} & \longrightarrow & \mathcal{P}\left(\left(\frac{\mathbb{F}_{q}[X]}{\langle X^{m} \rangle>}\right)^{n}\right) \\
 & (P_{1}(X)+\langle X^{m-1} \rangle,\ldots,P_{n}(X)+\langle X^{m-1} \rangle) & \longmapsto & \{(Q_{1}(X)+\langle X^{m} \rangle,\ldots,Q_{n}(X)+\langle X^{m} \rangle),~\exists c_{i} \in \mathbb{F}_{q} \\
 &                                                &             & Q_{i}(X)+\langle X^{m} \rangle=P_{i}(X)+c_{i}X^{m-1}+\langle X^{m} \rangle\}.  \\
\end{array}\]

\noindent We have the following equality:
\[\bigsqcup_{z \in R_{n}^{0}\left(\mathbb{F}_{q}[X]/\langle X^{m-1} \rangle\right)} \tilde{\iota}(z)= \bigsqcup_{c \in \mathbb{F}_{q}} R_{n}^{cX^{m-1}+\langle X^{m} \rangle}\left(\frac{\mathbb{F}_{q}[X]}{\langle X^{m} \rangle}\right).\]

\noindent Besides, $\frac{\mathbb{F}_{q}[X]}{\langle X^{m} \rangle}$ is a local ring and for all $c \in \mathbb{F}_{q}$ $cX^{m-1}+\langle X^{m} \rangle$ is non invertible (since $m \geq 2$). So, by proposition \ref{34}, we have for all $c \in \mathbb{F}_{q}$:
\[\left|R_{n}^{0}\left(\frac{\mathbb{F}_{q}[X]}{\langle X^{m} \rangle}\right)\right|=\left|R_{n}^{cX^{m-1}+\langle X^{m} \rangle}\left(\frac{\mathbb{F}_{q}[X]}{\langle X^{m} \rangle}\right)\right|.\] 

\noindent If we collect all the preceding elements, we have:
\begin{eqnarray*}
q~r_{n}\left(\frac{\mathbb{F}_{q}[X]}{\langle X^{m} \rangle}\right)=q\left|R_{n}^{0}\left(\frac{\mathbb{F}_{q}[X]}{\langle X^{m} \rangle}\right)\right| &=& \sum_{c \in \mathbb{F}_{q}} \left|R_{n}^{cX^{m-1}+\langle X^{m} \rangle}\left(\frac{\mathbb{F}_{q}[X]}{\langle X^{m} \rangle}\right)\right| \\
                                                                                                                                       &=& \left|\bigsqcup_{c \in \mathbb{F}_{q}} R_{n}^{cX^{m-1}+\langle X^{m} \rangle}\left(\frac{\mathbb{F}_{q}[X]}{\langle X^{m} \rangle}\right)\right| \\
																		                                                                                                   &=& \left|\bigsqcup_{z \in R_{n}^{0}\left(\frac{\mathbb{F}_{q}[X]}{\langle X^{m-1} \rangle}\right)} \tilde{\iota}(z)\right| \\
                                                                                                                                       &=& \sum_{z \in R_{n}^{0}\left(\frac{\mathbb{F}_{q}[X]}{\langle X^{m-1} \rangle}\right)} \left|\tilde{\iota}(z)\right| \\
                                                                                                                                       &=& \sum_{z \in R_{n}^{0}\left(\frac{\mathbb{F}_{q}[X]}{\langle X^{m-1} \rangle}\right)} q^{n} \\
																		                                                                                                   &=& q^{n}r_{n}\left(\frac{\mathbb{F}_{q}[X]}{\langle X^{m-1} \rangle}\right).
\end{eqnarray*}

\noindent Hence, $r_{n}\left(\frac{\mathbb{F}_{q}[X]}{\langle X^{m} \rangle}\right)=q^{n-1}r_{n}\left(\frac{\mathbb{F}_{q}[X]}{\langle X^{m-1} \rangle}\right)$.

\end{proof}

\begin{proof}[Second proof of corollary \ref{prin3}]

i) Let $A=\mathbb{F}_{q}$ and $n \geq 3$. $A$ is a finite local ring. First, we suppose $n$ odd. By theorem \ref{prin2}, $r_{n}(A)=q^{n}-(q-1)w_{n+2,A}^{1_{A}}$. By theorem \ref{26}, 
\[r_{n}(A)=q^{n}-(q-1)\frac{q^{n+1}-1}{q^{2}-1}=q^{n}-\frac{q^{n+1}-1}{q+1}=\frac{q^{n}+1}{q+1}.\]
\noindent Now, we suppose $n=2m$ even. By theorem \ref{prin2}, $r_{n}(A)=q^{n}-\sum_{u \in U(A)} \left|R_{n}^{u}(A)\right|$. We suppose $p \neq 2$. If we use lemma \ref{35}, and theorems \ref{prin2}, \ref{26} and \ref{27}, we obtain:
\begin{eqnarray*}
r_{n}(A) &=& q^{n}-\sum_{u \neq 0,1,-1} \frac{(q^{m+1}-1)(q^{m}-1)}{q^{2}-1}-q^{m}-2\frac{(q^{m+1}-1)(q^{m}-1)}{q^{2}-1} \\
        &=& q^{n}-(q-3)\frac{(q^{m+1}-1)(q^{m}-1)}{q^{2}-1}-q^{m}-2\frac{(q^{m+1}-1)(q^{m}-1)}{q^{2}-1} \\
				&=& q^{n}-q^{m}-\frac{(q^{m+1}-1)(q^{m}-1)(q-3+2)}{q^{2}-1} \\
				&=& q^{n}-q^{m}-\frac{(q^{m+1}-1)(q^{m}-1)}{q+1} \\
				&=& q^{n}-\frac{q^{n+1}+1}{q+1} \\
				&=& \frac{q^{n}-1}{q+1}. \\
\end{eqnarray*}

\noindent We suppose $p=2$. If we use lemma \ref{35}, and theorems \ref{prin2} and \ref{26}, we obtain:
\begin{eqnarray*}
r_{n}(A) &=& q^{n}-\sum_{u \neq 0,1} \frac{(q^{m+1}-1)(q^{m}-1)}{q^{2}-1}-q^{m}-\frac{(q^{m+1}-1)(q^{m}-1)}{q^{2}-1} \\
        &=& q^{n}-(q-2)\frac{(q^{m+1}-1)(q^{m}-1)}{q^{2}-1}-q^{m}-\frac{(q^{m+1}-1)(q^{m}-1)}{q^{2}-1} \\
				&=& q^{n}-q^{m}-\frac{(q^{m+1}-1)(q^{m}-1)(q-2+1)}{q^{2}-1} \\
				&=& q^{n}-q^{m}-\frac{(q^{m+1}-1)(q^{m}-1)}{q+1} \\
				&=& q^{n}-\frac{q^{n+1}+1}{q+1} \\
				&=& \frac{q^{n}-1}{q+1}.
\end{eqnarray*}

\noindent Hence, $r_{n}(A)=\frac{q^{n}+(-1)^{n+1}}{q+1}$. Besides, this formula is still true if $n=1$ or $n=2$.
\\
\\ii) Let $p$ be a prime number, $m \geq 2$ and $n \geq 1$, $A=\mathbb{Z}/p^{m}\mathbb{Z}$. By lemma \ref{36}, we have: 
\[r_{n}(A)=p^{n-1}r_{n}(\mathbb{Z}/p^{m-1}\mathbb{Z})=\ldots=p^{(m-1)(n-1)}r_{n}(\mathbb{Z}/p\mathbb{Z}).\]
\noindent By the case i), we have $r_{n}(A)=\frac{p^{(m-1)(n-1)}(p^{n}+(-1)^{n+1})}{p+1}$.
\\
\\iii) Let $q$ be a power of a prime number, $m \geq 2$ and $n \geq 1$, $A=\frac{\mathbb{F}_{q}[X]}{\langle X^{m} \rangle}$. By lemma \ref{37}, we have: 
\[r_{n}(A)=q^{n-1}r_{n}\left(\frac{\mathbb{F}_{q}[X]}{\langle X^{m-1} \rangle}\right)=\ldots=q^{(m-1)(n-1)}r_{n}\left(\frac{\mathbb{F}_{q}[X]}{\langle X \rangle}\right)=q^{(m-1)(n-1)}r_{n}(\mathbb{F}_{q}).\]
\noindent By the case i), we have $r_{n}(A)=\frac{q^{(m-1)(n-1)}(q^{n}+(-1)^{n+1})}{q+1}$.

\end{proof}

\begin{remark}
{\rm This result allows us to have the number of roots of the continuant over any finite commutative and unitary rings of order $p^{2}$, since these are isomorphic to one of the following rings: $\mathbb{F}_{p^{2}}$, $\mathbb{Z}/p^{2}\mathbb{Z}$, $\mathbb{Z}/p\mathbb{Z} \times \mathbb{Z}/p\mathbb{Z}$, $\frac{\mathbb{F}_{p}[X]}{\langle X^{2} \rangle}$ (see \cite{F}).
}
\end{remark}

\noindent Now, we can give another proof of the counting formulae for $w_{n+2}^{1}\left(\frac{\mathbb{F}_{q}[X]}{\langle X^{m} \rangle}\right)$ and $w_{n+2}^{\overline{1}}(\mathbb{Z}/p^{m}\mathbb{Z})$.

\begin{proof}[Second proof of corollary \ref{cor}]

i) Let $q$ be a power of a prime number, $m \geq 2$ and $n \geq 1$ odd, $A=\frac{\mathbb{F}_{q}[X]}{\langle X^{m} \rangle}$. $A$ is a local ring (proposition \ref{22}). By theorem \ref{prin2} and the second proof of corollary \ref{prin3} iii), we have the following equality:
\[\frac{q^{mn}-q^{m-1}(q-1)w_{n+2}^{1_{A}}(A)}{q^{m-1}}=r_{n}(A)=\frac{q^{(m-1)(n-1)}(q^{n}+1)}{q+1}.\]

\noindent Hence, we obtain the formula below:
\[w_{n+2}^{1_{A}}(A)=\frac{q^{nm}(q+1)-q^{(m-1)(n-1)+(m-1)}(q^{n}+1)}{q^{m-1}(q+1)(q-1)}=\frac{q^{nm+1}+q^{nm}-q^{nm}-q^{nm-n}}{q^{m-1}(q^{2}-1)}=\frac{q^{nm+1}-q^{n(m-1)}}{q^{m-1}(q^{2}-1)}.\]
\\
\\ii) Let $p$ be a prime number, $m \geq 2$ and $n \geq 1$ odd, $A=\mathbb{Z}/p^{m}\mathbb{Z}$. $A$ is a local ring (proposition \ref{22}). By theorem \ref{prin2} and the second proof of corollary \ref{prin3} ii), we have the following equality:
\[\frac{p^{mn}-p^{m-1}(p-1)w_{n+2}^{\overline{1}}(A)}{p^{m-1}}=r_{n}(A)=\frac{p^{(m-1)(n-1)}(p^{n}+1)}{p+1}.\]

\noindent Hence, we obtain the formula below:
\[w_{n+2}^{\overline{1}}(A)=\frac{p^{nm}(p+1)-p^{(m-1)(n-1)+(m-1)}(p^{n}+1)}{p^{m-1}(p+1)(p-1)}=\frac{p^{nm+1}+p^{nm}-p^{nm}-p^{nm-n}}{p^{m-1}(p^{2}-1)}=\frac{p^{nm+1}-p^{n(m-1)}}{p^{m-1}(p^{2}-1)}.\]

\end{proof}

\subsection{Numerical applications}
\label{chap34}

\noindent Using the Chinese remainder theorem, we can easily prove the following result:

\begin{proposition}

Let $N=p_{1}^{\alpha_{1}}\ldots p_{r}^{\alpha_{r}}$ with $p_{i}$ distinct prime numbers and $\alpha_{i} \geq 1$. Let $n \geq 1$.  
\[r_{n}(\mathbb{Z}/N\mathbb{Z})=\prod_{i=1}^{r} r_{n}(\mathbb{Z}/p_{i}^{\alpha_{i}}\mathbb{Z}).\]

\end{proposition}

\noindent The tables below give some values of $r_{n}(A)$.

\hfill\break

\begin{center}
\begin{tabular}{|c|c|c|c|c|c|c|c|c|c|c|}
\hline
  \multicolumn{1}{|c|}{\backslashbox{$n$}{\vrule width 0pt height 1.25em$A$}}     & $\mathbb{F}_{2}$ & $\mathbb{F}_{3}$ & $\mathbb{F}_{4}$ & $\mathbb{Z}/4\mathbb{Z}$ & $\mathbb{F}_{5}$ & $\mathbb{Z}/6\mathbb{Z}$ & $\mathbb{F}_{7}$ & $\mathbb{Z}/8\mathbb{Z}$ & $\mathbb{F}_{9}$ & $\mathbb{Z}/9\mathbb{Z}$   \rule[-7pt]{0pt}{18pt} \\
  \hline
  3 & 3 & 7 & 13 & 12 & 21 & 21 & 43 & 48 & 73 & 63   \rule[-7pt]{0pt}{18pt} \\
	\hline
  4 & 5 & 20 & 51 & 40 & 104 & 100 & 300 & 320 & 656 & 540 \rule[-7pt]{0pt}{18pt} \\
	\hline
  5  & 11 & 61 & 205 & 176 & 521 & 671 & 2101 & 2816 & 5905 & 4941 \rule[-7pt]{0pt}{18pt} \\
	\hline
  6  & 21 & 182 & 819 & 672 & 2604 & 3822 & 14~706 & 21~504 & 53~144 & 44~226 \rule[-7pt]{0pt}{18pt} \\
  \hline
	7  & 43 & 547 & 3277 & 2752 & 13~021 & 23~521 & 102~943 & 176~128 & 478~297 & 398~763  \rule[-7pt]{0pt}{18pt} \\
	\hline
	8  & 85 & 1640 & 13~107 & 10~880 & 65~104 & 139~400 & 720~600 & 1~392~640 & 4~304~672 & 3~586~680 \rule[-7pt]{0pt}{18pt} \\
	\hline
	
\end{tabular}
\end{center}

\hfill\break

\begin{center}
\begin{tabular}{|c|c|c|c|c|c|c|}
\hline
  \multicolumn{1}{|c|}{\backslashbox{$n$}{\vrule width 0pt height 1.25em$A$}}    & $\mathbb{Z}/10\mathbb{Z}$ & $\mathbb{F}_{11}$ & $\mathbb{Z}/12\mathbb{Z}$ & $\mathbb{F}_{16}$ & $\mathbb{Z}/16\mathbb{Z}$ & $\frac{\mathbb{F}_{4}[X]}{\langle X^{2} \rangle}$   \rule[-7pt]{0pt}{18pt} \\
  \hline
  3  & 63 & 111 & 84  & 241 & 192 & 208 \rule[-7pt]{0pt}{18pt} \\
	\hline
  4  & 520 & 1220 & 800  & 3855 & 2560 & 3264 \rule[-7pt]{0pt}{18pt} \\
	\hline
  5  & 5731 & 13~421 & 10~736 & 61~681 & 45~056 & 52~480 \rule[-7pt]{0pt}{18pt} \\
	\hline
  6  & 54~684 & 147~630 & 122~304 & 986~895 & 688~128 & 838~656  \rule[-7pt]{0pt}{18pt} \\
  \hline
	7  & 559~903 & 1~623~931 & 1~505~344 & 15~790~321 & 11~272~192 & 13~422~592 \rule[-7pt]{0pt}{18pt} \\
	\hline
	8  & 5~533~840 & 17~863~240 & 17~843~200 & 252~645~135 & 178~257~920 & 214~745~088 \rule[-7pt]{0pt}{18pt} \\
	\hline
	
\end{tabular}
\end{center}

\hfill\break

\newcolumntype{N}{@{}m{0pt}@{}}
\begin{center}
\begin{tabular}{|>{\centering\arraybackslash}m{2.65cm}|>{\centering\arraybackslash}m{1cm}|>{\centering\arraybackslash}m{1cm}|>{\centering\arraybackslash}m{1cm}|>{\centering\arraybackslash}m{1cm}|>{\centering\arraybackslash}m{1.5cm}|>{\centering\arraybackslash}m{1.8cm}|>{\centering\arraybackslash}m{1.8cm}|N|}
\hline
  \diagbox[width=3cm]{$A$}{$n$} & 2 & 3 & 4 & 5 & 6 & 7 & 8 &  \rule[-7pt]{0pt}{18pt} \\
		\hline
  $\displaystyle{\frac{(\mathbb{Z}/4\mathbb{Z})[X]}{\langle X^{2},\overline{2}X \rangle}}$ & 4 & 48 & 320 & 2816 & 21~504 & 176~128 & 1~392~640 & \rule[-7pt]{0pt}{28pt} \\
  \hline
  $\displaystyle{\frac{(\mathbb{Z}/2\mathbb{Z})[X,Y]}{\langle X^{2},Y^{2} \rangle}}$ & 8 & 192 & 2560 & 45~056 & 688~128 & 11~272~192 & 178~257~920 & \rule[-7pt]{0pt}{28pt} \\
	\hline
  $\displaystyle{\frac{(\mathbb{Z}/4\mathbb{Z})[X]}{\langle X^{2}+X+\overline{1} \rangle}}$ & 12 & 208 & 3264 & 52~480 & 838~656 & 13~422~592 & 214~745~088 & \rule[-7pt]{0pt}{28pt} \\
	\hline
	
\end{tabular}
\end{center}

\hfill\break

\section{Number of $\lambda$-quiddities over $\mathbb{Z}/p^{m}\mathbb{Z}$}

The aim of this section is to prove theorem \ref{prin4}. Note that the resolution of equation \eqref{p} over $\mathbb{Z}/N\mathbb{Z}$ is closely linked to the different writings of the elements of the following congruence subgroup: \[\hat{\Gamma}(N):=\{C \in SL_{2}(\mathbb{Z}), C \equiv \pm Id [N]\}.\]
\noindent Indeed, we know that all the matrices of the modular group can be written in the form $M_{n}(a_{1},\ldots,a_{n})$, with $a_{i}$ a positive integer. Since this representation is not unique, we are naturally led to investigate all representations of this form for a given matrix, or more generally, for a set of matrices.

\subsection{Preliminary lemmas}

\begin{lemma}
\label{41}

Let $n, m \geq  2$, $p$ be a prime number, and $N=p^{m}$. Let $\epsilon \in \{-1,1\}$. 
\begin{itemize}
\item $n=2l \geq 4$ even.
\[w_{n+2}^{\overline{\epsilon}}(\mathbb{Z}/N\mathbb{Z})=p^{n-1}w_{n+2}^{\epsilon}(\mathbb{Z}/p^{m-1}\mathbb{Z})+p^{m(l-1)}\left|R_{2}^{(-1)^{l}\overline{\epsilon}}(\mathbb{Z}/N \mathbb{Z})\right|-p^{m(l-1)-1}\sum_{k=0}^{p-1} \left|R_{2}^{(-1)^{l}\overline{(\epsilon+kp^{m-1})}^{(-1)^{l+1}}}(\mathbb{Z}/N \mathbb{Z})\right|.\]
\item $n \geq 3$ odd.
\[w_{n+2}^{\overline{\epsilon}}(\mathbb{Z}/N \mathbb{Z})=p^{n-1}w_{n+2}^{\epsilon}(\mathbb{Z}/p^{m-1}\mathbb{Z}).\]
\end{itemize}

\end{lemma}

\begin{proof}

Let $N:=p^{m}$, $U:=U(\mathbb{Z}/N\mathbb{Z})$, $u \in \mathbb{Z}$ such that $\overline{u} \in U$ and $W:=(\mathbb{Z}/N\mathbb{Z}-U)-\{\overline{0}\}$. Let $(a,b) \in (\mathbb{Z}/p^{m}\mathbb{Z})^{2}$, we set $S_{n}^{a,b}(N):=\{(a_{1},\ldots,a_{n}) \in (\mathbb{Z}/p^{m}\mathbb{Z})^{n}, K_{n}(a_{1},\ldots,a_{n})=a~{\rm and}~K_{n-1}(a_{1},\ldots,a_{n-1})=b\}$ and, if $a$ is invertible, $B_{a}=\begin{pmatrix}
                 a & \overline{0} \\
								\overline{0} & a^{-1}
								\end{pmatrix}$. To simplify the notations, we will write $R_{n}^{-u}(\mathbb{Z}/p^{m-1}\mathbb{Z})$ instead of $R_{n}^{-u+p^{m-1}\mathbb{Z}}(\mathbb{Z}/p^{m-1}\mathbb{Z})$ and $w_{n+2}^{u}(\mathbb{Z}/p^{m-1}\mathbb{Z})$ instead of $w_{n+2}^{u+p^{m-1}\mathbb{Z}}(\mathbb{Z}/p^{m-1}\mathbb{Z})$.
\\
\\We define the following application:
\[\begin{array}{ccccc} 
\iota : & (\mathbb{Z}/p^{m-1}\mathbb{Z})^{n} & \longrightarrow & \mathcal{P}((\mathbb{Z}/p^{m}\mathbb{Z})^{n}) \\
 & (a_{1}+\frac{N}{p}\mathbb{Z},\ldots,a_{n}+\frac{N}{p}\mathbb{Z}) & \longmapsto & \{(\overline{b_{1}},\ldots,\overline{b_{n}}),\exists k_{i} \in [\![0;p-1]\!]~\overline{b_{i}}=\overline{a_{i}+\frac{k_{i}N}{p}}\}  \\
\end{array}.\]

\noindent We have $\bigsqcup_{z \in R_{n}^{-u}(\mathbb{Z}/p^{m-1}\mathbb{Z})} \iota(z)= \bigsqcup_{k=0}^{p-1} R_{n}^{\overline{-u-\frac{kN}{p}}}(\mathbb{Z}/N \mathbb{Z})$ (the proof is the same as that carried out in lemma \ref{36}).
\\
\\Let $k \in [\![0;p-1]\!]$ and $x=\overline{u+\frac{kN}{p}}$. We have the following equality:

\[R_{n}^{-x}(\mathbb{Z}/N \mathbb{Z})=\bigsqcup_{v \in U} S_{n}^{-x,v}(N)  \bigsqcup_{w \in W} S_{n}^{-x,w}(N) \bigsqcup S_{n}^{-x,\overline{0}}(N).\] 

\noindent Let $v \in U$. Suppose $S_{n}^{-\overline{u},v}(N) \neq \emptyset$. The following application is well-defined (the proof is the same as that carried out in proposition \ref{34}):
\[\begin{array}{ccccc} 
\sigma_{n,u,v} : & S_{n}^{-\overline{u},v}(N) & \longrightarrow & S_{n}^{-x,v}(N) \\
  & (a_{1},\ldots,a_{n}) & \longmapsto & (a_{1},\ldots,a_{n-1},a_{n}-v^{-1}(\overline{kp^{m-1}}))  \\
\end{array}.\]
\noindent Besides, $S_{n}^{-x,v}(N) \neq \emptyset$ and we can define :
\[\begin{array}{ccccc} 
\tilde{\sigma}_{n,u,v} : & S_{n}^{-x,v}(N) & \longrightarrow & S_{n}^{-\overline{u},v}(N) \\
 & (a_{1},\ldots,a_{n}) & \longmapsto &  (a_{1},\ldots,a_{n-1},a_{n}+v^{-1}(\overline{kp^{m-1}}))  \\
\end{array}.\]

\noindent $\sigma_{n,u,v}$ and $\tilde{\sigma}_{n,u,v}$ are reciprocal bijections. Hence, $\left|S_{n}^{-\overline{u},v}(N)\right|=\left|S_{n}^{-x,v}(N)\right|$.
\\
\\\noindent Let $w \in W$. There exists $h \in \mathbb{Z}$ and $l \in [\![1;m-1]\!]$ such that $w=\overline{hp^{l}}$ and $p$ does not divide $h$. Suppose $S_{n}^{-\overline{u},w}(N) \neq \emptyset$. The following application is well-defined (the proof is the same as that carried out in proposition \ref{34}):
\[\begin{array}{ccccc} 
\tau_{n,u,w} : & S_{n}^{-\overline{u},w}(N) & \longrightarrow & S_{n}^{-x,w}(N) \\
  & (a_{1},\ldots,a_{n}) & \longmapsto & (a_{1},\ldots,a_{n-1},a_{n}-(\overline{h^{-1}kp^{m-1-l}}))  \\
\end{array}.\]
\noindent Besides, $S_{n}^{-x,w}(N) \neq \emptyset$ and we can define:
\[\begin{array}{ccccc} 
\tilde{\tau}_{n,u,w} : & S_{n}^{-x,w}(N) & \longrightarrow & S_{n}^{-\overline{u},w}(N) \\
 & (a_{1},\ldots,a_{n}) & \longmapsto &  (a_{1},\ldots,a_{n-1},a_{n}+(\overline{h^{-1}kp^{m-1-l}}))  \\
\end{array}.\]

\noindent $\tau_{n,u,w}$ and $\tilde{\tau}_{n,u,w}$ are reciprocal bijections. Thus, $\left|S_{n}^{-\overline{u},w}(N)\right|=\left|S_{n}^{-x,w}(N)\right|$.
\\
\\Moreover, let $z \in U$, we have the following equality:
\[S_{n}^{z,\overline{0}}=\bigsqcup_{y \in \mathbb{Z}/N\mathbb{Z}} \underbrace{\left\{(a_{1},\ldots,a_{n}) \in (\mathbb{Z}/N\mathbb{Z})^{n},~M_{n}(a_{1},\ldots,a_{n})=\begin{pmatrix}
                 z & y \\
								\overline{0} & z^{-1}
								\end{pmatrix}\right\}}_{T_{n,z,y}},\]
since ${\rm det}(M_{n}(a_{1},\ldots,a_{n}))=\overline{1}$. Let $y \in \mathbb{Z}/N\mathbb{Z}$. Suppose $\Omega_{n}^{B_{z}}(\mathbb{Z}/N\mathbb{Z}) \neq \emptyset$. The following application is well-defined:
\[\begin{array}{ccccc} 
\zeta_{n,z,y} : & \Omega_{n}^{B_{z}}(\mathbb{Z}/N\mathbb{Z}) & \longrightarrow & T_{n,z,y}(N) \\
  & (a_{1},\ldots,a_{n}) & \longmapsto & (a_{1},\ldots,a_{n-1},a_{n}+yz)  \\
\end{array}.\]
\noindent Besides, $T_{n,z,y}(N) \neq \emptyset$ and we can define:
\[\begin{array}{ccccc} 
\tilde{\zeta}_{n,z,y} : & T_{n,z,y}(N) & \longrightarrow & \Omega_{n}^{B_{z}}(\mathbb{Z}/N\mathbb{Z}) \\
 & (a_{1},\ldots,a_{n}) & \longmapsto &  (a_{1},\ldots,a_{n-1},a_{n}-yz)  \\
\end{array}.\]

\noindent $\zeta_{n,z,y}$ and $\tilde{\zeta}_{n,z,y}$ are reciprocal bijections. Hence, $\left|\Omega_{n}^{B_{z}}(\mathbb{Z}/N\mathbb{Z})\right|=\left|T_{n,z,y}(N)\right|$. So, if $n \geq 3$, we have, by lemma \ref{32}:
\[\left|S_{n}^{-\overline{u},\overline{0}}\right|=p^{m}\left|\Omega_{n}^{B_{-\overline{u}}}(\mathbb{Z}/N\mathbb{Z})\right|=p^{m}\left|R_{n-2}^{\overline{u}^{-1}}(\mathbb{Z}/N\mathbb{Z})\right|.\]

\noindent If we collect all the preceding elements and if we use lemma \ref{32}, we have for all $n \geq 3$:
\begin{eqnarray*}
\left|R_{n}^{-x}(\mathbb{Z}/N \mathbb{Z})\right| &=& \sum_{v \in U} \left|S_{n}^{-x,v}(N)\right| +  \sum_{w \in W} \left|S_{n}^{-x,w}(N)\right| + \left|S_{n-2}^{-x,\overline{0}}(N)\right| \\
                                                 &=& \sum_{v \in U} \left|S_{n}^{-\overline{u},v}(N)\right| +  \sum_{w \in W} \left|S_{n}^{-\overline{u},w}(N)\right| + \sum_{y \in \mathbb{Z}/N\mathbb{Z}} \left|T_{n,-x,y}(N)\right| \\
																		             &=& \sum_{v \in U} \left|S_{n}^{-\overline{u},v}(N)\right| +  \sum_{w \in W} \left|S_{n}^{-\overline{u},w}(N)\right| + \sum_{y \in \mathbb{Z}/N\mathbb{Z}} \left|\Omega_{n}^{B_{-x}}(\mathbb{Z}/N \mathbb{Z})\right| \\
																								 &=& \left|\bigsqcup_{v \in U} S_{n}^{-\overline{u},v}(N) \bigsqcup_{w \in W} S_{n}^{-\overline{u},w}(N) \bigsqcup S_{n}^{-\overline{u},\overline{0}}(N)\right|-\left|S_{n}^{-\overline{u},\overline{0}}(N)\right| + p^{m} \left|\Omega_{n}^{B_{-x}}(\mathbb{Z}/N \mathbb{Z})\right| \\
																								 &=& \left|R_{n}^{-\overline{u}}(\mathbb{Z}/N \mathbb{Z})\right|+p^{m}\left(\left|R_{n-2}^{x^{-1}}(\mathbb{Z}/N \mathbb{Z})\right|-\left|R_{n-2}^{\overline{u}^{-1}}(\mathbb{Z}/N \mathbb{Z})\right|\right).																					 
\end{eqnarray*}

\noindent Now, we suppose $\overline{u}=\overline{\epsilon}$. We have $\left|R_{n}^{-x}(\mathbb{Z}/N \mathbb{Z})\right|=\left|R_{n}^{-\overline{\epsilon}}(\mathbb{Z}/N \mathbb{Z})\right|+p^{m}\left(\left|R_{n-2}^{x^{-1}}(\mathbb{Z}/N \mathbb{Z})\right|-\left|R_{n-2}^{\overline{\epsilon}}(\mathbb{Z}/N \mathbb{Z})\right|\right)$. Since, $x^{-1}=(\overline{\epsilon+kp^{m-1}})^{-1}=\overline{\epsilon+(p-k)p^{m-1}}$, the previous recursive formula is still true if we replace $-x$ by $x^{-1}$, that is to say $\left|R_{n-2}^{-(-x^{-1})}(\mathbb{Z}/N \mathbb{Z})\right|=\left|R_{n-2}^{\overline{\epsilon}}(\mathbb{Z}/N \mathbb{Z})\right|+p^{m}\left(\left|R_{n-4}^{-x}(\mathbb{Z}/N \mathbb{Z})\right|-\left|R_{n-4}^{-\overline{\epsilon}}(\mathbb{Z}/N \mathbb{Z})\right|\right)$ (if $n \geq 5$).
\\
\\Suppose $n=2l \geq 4$ even. We have the following equalities:
\[\left|R_{n}^{-x}(\mathbb{Z}/N \mathbb{Z})\right|=\ldots=\left|R_{n}^{-\overline{\epsilon}}(\mathbb{Z}/N \mathbb{Z})\right|+p^{m(l-1)}\left(\left|R_{2}^{(-1)^{l}x^{(-1)^{l+1}}}(\mathbb{Z}/N \mathbb{Z})\right|-\left|R_{2}^{(-1)^{l}\overline{\epsilon}}(\mathbb{Z}/N \mathbb{Z})\right|\right).\]

\noindent Hence, 
\begin{eqnarray*} 
\sum_{k=0}^{p-1} \left|R_{n}^{\overline{-\epsilon-kp^{m-1}}}(\mathbb{Z}/N \mathbb{Z})\right| &=& p\left|R_{n}^{-\overline{\epsilon}}(\mathbb{Z}/N \mathbb{Z})\right|-p^{m(l-1)+1}\left|R_{2}^{(-1)^{l}\overline{\epsilon}}(\mathbb{Z}/N \mathbb{Z})\right| \\
                                                                                      & & +p^{m(l-1)}\sum_{k=0}^{p-1} \left|R_{2}^{(-1)^{l}\overline{(\epsilon+kp^{m-1}})^{(-1)^{l+1}}}(\mathbb{Z}/N \mathbb{Z})\right|.
\end{eqnarray*}																																											

\noindent Besides, 
\[\sum_{k=0}^{p-1} \left|R_{n}^{\overline{-\epsilon-kp^{m-1}}}(\mathbb{Z}/N \mathbb{Z})\right|=\left|\bigsqcup_{z \in R_{n}^{-\epsilon}(\mathbb{Z}/p^{m-1}\mathbb{Z})} \iota(z)\right|=\sum_{z \in R_{n}^{-\epsilon}(\mathbb{Z}/p^{m-1}\mathbb{Z})} \left|\iota(z)\right|=p^{n}\left|R_{n}^{-\epsilon}(\mathbb{Z}/p^{m-1}\mathbb{Z})\right|.\]

\noindent Thus, by lemma \ref{32}, we have:
\begin{eqnarray*}
w_{n+2}^{\overline{\epsilon}}(\mathbb{Z}/N \mathbb{Z}) &=& \left|R_{n}^{-\overline{\epsilon}}(\mathbb{Z}/N \mathbb{Z})\right| \\
                                                       &=& p^{n-1}\left|R_{n}^{-\epsilon}(\mathbb{Z}/p^{m-1}\mathbb{Z})\right|+p^{m(l-1)}\left|R_{2}^{(-1)^{l}\overline{\epsilon}}(\mathbb{Z}/N \mathbb{Z})\right| \\
																								       & &-p^{m(l-1)-1}\sum_{k=0}^{p-1} \left|R_{2}^{(-1)^{l}\overline{(\epsilon+kp^{m-1}})^{(-1)^{l+1}}}(\mathbb{Z}/N \mathbb{Z})\right| \\
																								       &=& p^{n-1}w_{n+2}^{\epsilon}(\mathbb{Z}/p^{m-1}\mathbb{Z})+p^{m(l-1)}\left|R_{2}^{(-1)^{l}\overline{\epsilon}}(\mathbb{Z}/N \mathbb{Z})\right| \\
																								       & & -p^{m(l-1)-1}\sum_{k=0}^{p-1} \left|R_{2}^{(-1)^{l}\overline{(\epsilon+kp^{m-1}})^{(-1)^{l+1}}}(\mathbb{Z}/N \mathbb{Z})\right|.
\end{eqnarray*}

\noindent Suppose $n=2l+1 \geq 3$ odd. Since, $\left|R_{1}^{a}(\mathbb{Z}/N \mathbb{Z})\right|=1$ for all $a \in \mathbb{Z}/N \mathbb{Z}$, the same elements give us $\left|R_{n}^{-x}(\mathbb{Z}/N \mathbb{Z})\right|=\left|R_{n}^{-\overline{\epsilon}}(\mathbb{Z}/N \mathbb{Z})\right|$ and $w_{n+2}^{\overline{\epsilon}}(\mathbb{Z}/N \mathbb{Z})=p^{n-1}w_{n+2}^{\epsilon}(\mathbb{Z}/p^{m-1}\mathbb{Z})$.						

\end{proof}

\begin{lemma}
\label{42}

Let $m \geq 2$, $p$ be a prime number, and $k \in [\![0;p-1]\!]$. 
\begin{itemize}
\item $p$ odd. $w_{4}^{-\overline{1}}(\mathbb{Z}/p^{m}\mathbb{Z})=\left|R_{2}^{\overline{1+kp^{m-1}}}(\mathbb{Z}/p^{m}\mathbb{Z})\right|=p^{m-1}(p-1)$.
\item $\sum_{j=0}^{p-1} \left|R_{2}^{\overline{-1-jp^{m-1}}}(\mathbb{Z}/p^{m}\mathbb{Z})\right|=p^{m}(m(p-1)+1)$; 
\\$w_{4}^{\overline{1}}(\mathbb{Z}/p^{m}\mathbb{Z})=\left|R_{2}^{\overline{-1}}(\mathbb{Z}/p^{m}\mathbb{Z})\right|=p^{m-1}(p(m+1)-m)$.
\item $m \geq 3$, $w_{4}^{-\overline{1}}(\mathbb{Z}/2^{m}\mathbb{Z})=\left|R_{2}^{\overline{1}}(\mathbb{Z}/2^{m}\mathbb{Z})\right|=\left|R_{2}^{\overline{1+2^{m-1}}}(\mathbb{Z}/2^{m}\mathbb{Z})\right|=2^{m}$.
\end{itemize}

\end{lemma}

\begin{proof}

i) $x=\overline{2+kp^{m-1}}$ is invertible (since $p$ is odd). Hence, $R_{2}(\overline{1+kp^{m-1}})=\{(a,a^{-1}x), a \in U(\mathbb{Z}/p^{m}\mathbb{Z})\}$. Thus, $\left|R_{2}^{\overline{1+kp^{m-1}}}(\mathbb{Z}/p^{m}\mathbb{Z})\right|=\left|U(\mathbb{Z}/p^{m}\mathbb{Z})\right|=p^{m-1}(p-1)$ and $w_{4}^{-\overline{1}}(\mathbb{Z}/p^{m}\mathbb{Z})=\left|R_{2}^{\overline{1}}(\mathbb{Z}/p^{m}\mathbb{Z})\right|$ (lemma \ref{32}).
\\
\\ii) $\bigsqcup_{j=0}^{p-1} R_{2}^{\overline{-1-jp^{m-1}}}(\mathbb{Z}/p^{m}\mathbb{Z})=\{(0+p^{m}\mathbb{Z},b+p^{m}\mathbb{Z}), 0\leq b<p^{m}\} \sqcup \{(a+p^{m}\mathbb{Z},0+p^{m}\mathbb{Z}), 0<a<p^{m}\} \sqcup\{(ap^{u}+p^{m}\mathbb{Z},bp^{v}+p^{m}\mathbb{Z}), 0 \leq u \leq m-1, m-1-u \leq v \leq m-1, p \nmid a, p \nmid b, 0<ap^{u}<p^{m}, 0<bp^{v}<p^{m}\}$.
\begin{eqnarray*}
\sum_{j=0}^{p-1} \left|R_{2}^{\overline{-1-jp^{m-1}}}(\mathbb{Z}/p^{m}\mathbb{Z})\right| &=& \left|\bigsqcup_{j=0}^{p-1} R_{2}^{\overline{-1-jp^{m-1}}}(\mathbb{Z}/p^{m}\mathbb{Z})\right| \\
	                                                          &=& p^{m}+p^{m}-1+\sum_{u=0}^{m-1} \sum_{v=m-1-u}^{m-1} p^{m-u-1}(p-1)p^{m-v-1}(p-1) \\
																														&=& 2p^{m}-1+(p-1)^{2}p^{2m-2}\sum_{u=0}^{m-1} p^{-u} \sum_{v=m-u-1}^{m-1} p^{-v}
\end{eqnarray*}	
\begin{eqnarray*}																													
\sum_{j=1}^{p-1} \left|R_{2}^{\overline{-1-jp^{m-1}}}(\mathbb{Z}/p^{m}\mathbb{Z})\right| &=& 2p^{m}-1+(p-1)^{2}p^{2m-2}\sum_{u=0}^{m-1} p^{-u} p^{u+1-m} \frac{1-p^{-u-1}}{1-p^{-1}} \\
	                                                          &=& 2p^{m}-1+(p-1)^{2}p^{m-1}\sum_{u=0}^{m-1} \frac{p-p^{-u}}{p-1} \\
	                                                          &=& 2p^{m}-1+(p-1)p^{m-1}\sum_{u=0}^{m-1} (p-p^{-u}) \\
	                                                          &=& 2p^{m}-1+(p-1)p^{m-1}\left(mp-\frac{1-p^{-m}}{1-p^{-1}}\right) \\
	                                                          &=& 2p^{m}-1+p^{m-1}(mp(p-1)-p+p^{-m+1}) \\
																														&=& p^{m}(m(p-1)+1).
\end{eqnarray*}

\noindent A similar calculation gives us $w_{4}^{\overline{1}}(\mathbb{Z}/p^{m}\mathbb{Z})=\left|R_{2}^{\overline{-1}}(\mathbb{Z}/p^{m}\mathbb{Z})\right|=p^{m-1}(p(m+1)-m)$.
\\
\\iii) Let $\epsilon \in \{0,1\}$. $R_{2}^{\overline{1+\epsilon 2^{m-1}}}(\mathbb{Z}/2^{m}\mathbb{Z})=\{(\overline{(2+\epsilon 2^{m-1})}a^{-1},a),~a \in U(\mathbb{Z}/2^{m}\mathbb{Z})\} \sqcup \{(a,(\overline{2+\epsilon 2^{m-1}})a^{-1}),~a \in U(\mathbb{Z}/2^{m}\mathbb{Z})\}$. Indeed, if $(\overline{a},\overline{b}) \in R_{2}^{\overline{1+\epsilon 2^{m-1}}}(\mathbb{Z}/2^{m}\mathbb{Z})$ then $ab$ is even but is not divisible by 4 (since $m \geq 3$). Hence, $a$ or $b$ is odd. Thus, $w_{4}^{-\overline{1}}(\mathbb{Z}/2^{m}\mathbb{Z})=\left|R_{2}^{\overline{1+\epsilon 2^{m-1}}}(\mathbb{Z}/2^{m}\mathbb{Z})\right|=2\times 2^{m-1}=2^{m}$.

\end{proof}

\subsection{Proof of the counting formulae}

If we collect all the preceding elements, we can prove the main result of this section.

\begin{proof}[Proof of theorem \ref{prin4} and third proof of the second point of corollary \ref{cor} ]

i) Let $n \geq 3$ odd and $\epsilon \in \{-1,1\}$. By lemma \ref{41} and theorem \ref{26}, we have the equalities below.
\[w_{n+2}^{\overline{\epsilon}}(\mathbb{Z}/N \mathbb{Z})=p^{n-1}w_{n+2}^{\epsilon}(\mathbb{Z}/p^{m-1}\mathbb{Z})=p^{(m-1)(n-1)}w_{n+2}^{\epsilon}(\mathbb{F}_{p})=\frac{p^{nm+1}-p^{n(m-1)}}{p^{m-1}(p^{2}-1)}.\] 

\noindent ii) Let $n=2l \geq 4$, $l$ even, $m \geq 2$ and $N=p^{m}$ with $p$ an odd prime number. By lemma \ref{41} and lemma \ref{42}, we have:
\begin{eqnarray*}
w_{n+2}^{\overline{1}}(\mathbb{Z}/p^{m}\mathbb{Z}) &=& p^{n-1}w_{n+2}^{1}(\mathbb{Z}/p^{m-1}\mathbb{Z})+p^{m(l-1)}\left|R_{2}^{\overline{1}}(\mathbb{Z}/N\mathbb{Z})\right|-p^{m(l-1)-1}\sum_{k=0}^{p-1} \left|R_{2}^{\overline{(1+kp^{m-1})}^{-1}}(\mathbb{Z}/N \mathbb{Z})\right| \\
                                                   &=& p^{n-1}w_{n+2}^{1}(\mathbb{Z}/p^{m-1}\mathbb{Z})+p^{m(l-1)}\left|R_{2}^{\overline{1}}(\mathbb{Z}/N\mathbb{Z})\right|-p^{m(l-1)-1}\sum_{k=0}^{p-1} \left|R_{2}^{\overline{(1+(p-k)p^{m-1})}}(\mathbb{Z}/N \mathbb{Z})\right| \\
                                                   &=& p^{n-1}w_{n+2}^{1}(\mathbb{Z}/p^{m-1}\mathbb{Z})+p^{m(l-1)}\left|R_{2}^{\overline{1}}(\mathbb{Z}/N\mathbb{Z})\right|-p^{m(l-1)-1}\sum_{k=0}^{p-1} \left|R_{2}^{\overline{1}}(\mathbb{Z}/N \mathbb{Z})\right| \\
																									 &=& p^{n-1}w_{n+2}^{1}(\mathbb{Z}/p^{m-1}\mathbb{Z}) \\
																									 &=& \ldots~= p^{(m-1)(n-1)} w_{2(l+1)}^{1}(\mathbb{Z}/p\mathbb{Z}). \\
\end{eqnarray*}

\noindent iii) Let $n=2l \geq 4$, $l$ even, $m \geq 2$ and $N=p^{m}$, with $p$ an odd prime number. By lemma \ref{41} and lemma \ref{42}, we have:
\begin{eqnarray*}
w_{n+2}^{\overline{-1}}(\mathbb{Z}/p^{m}\mathbb{Z}) &=& p^{n-1}w_{n+2}^{-1}(\mathbb{Z}/p^{m-1}\mathbb{Z})+p^{m(l-1)}\left|R_{2}^{\overline{-1}}(\mathbb{Z}/N\mathbb{Z})\right|-p^{m(l-1)-1}\sum_{k=0}^{p-1} \left|R_{2}^{\overline{-1+(p-k)p^{m-1}}}(\mathbb{Z}/N \mathbb{Z})\right|
\end{eqnarray*}
 \begin{eqnarray*}                                                   
w_{n+2}^{\overline{-1}}(\mathbb{Z}/p^{m}\mathbb{Z}) &=& p^{n-1}w_{n+2}^{-1}(\mathbb{Z}/p^{m-1}\mathbb{Z})+p^{m(l-1)}p^{m-1}(p(m+1)-m)-p^{m(l-1)-1}p^{m}(m(p-1)+1) \\
																									  &=& p^{n-1}w_{n+2}^{-1}(\mathbb{Z}/p^{m-1}\mathbb{Z})+p^{ml-1}((p(m+1)-m)-(m(p-1)+1)) \\
																									  &=& p^{n-1}w_{n+2}^{-1}(\mathbb{Z}/p^{m-1}\mathbb{Z})+p^{ml-1}(p-1) \\
																									  &=& \ldots~= p^{(m-1)(n-1)}w_{2(l+1)}^{-1}(\mathbb{Z}/p\mathbb{Z})+\sum_{j=0}^{m-2} p^{j(n-1)}p^{(m-j)l-1}(p-1) \\
																										&=& \ldots~= p^{(m-1)(n-1)}w_{2(l+1)}^{-1}(\mathbb{Z}/p\mathbb{Z})+p^{ml-1}(p-1)\sum_{j=0}^{m-2} p^{j(n-1-l)} \\
																									  &=& p^{(m-1)(n-1)}w_{2(l+1)}^{-1}(\mathbb{Z}/p\mathbb{Z})+\frac{p^{(m-1)(n-1-l)}-1}{p^{n-1-l}-1}p^{ml-1}(p-1).
\end{eqnarray*}

\noindent iv) Let $n=2l \geq 4$, $l$ even and $m \geq 3$. By lemma \ref{41} and lemma \ref{42}, we have:
\begin{eqnarray*}
w_{n+2}^{\overline{1}}(\mathbb{Z}/2^{m}\mathbb{Z}) &=& 2^{n-1}w_{n+2}^{1}(\mathbb{Z}/2^{m-1}\mathbb{Z})+2^{m(l-1)}\left|R_{2}^{\overline{1}}(\mathbb{Z}/2^{m}\mathbb{Z})\right|-2^{m(l-1)-1}\sum_{k=0}^{1} \left|R_{2}^{\overline{1+(2-k)2^{m-1}}}(\mathbb{Z}/2^{m}\mathbb{Z})\right| \\
                                                   &=& 2^{n-1}w_{n+2}^{1}(\mathbb{Z}/2^{m-1}\mathbb{Z})+2^{m(l-1)}2^{m}-2^{m(l-1)-1}\sum_{k=0}^{p-1} 2^{m} \\
																									 &=& 2^{n-1}w_{n+2}^{1}(\mathbb{Z}/2^{m-1}\mathbb{Z}) \\
																									 &=& \ldots~=2^{(m-2)(n-1)} w_{2(l+1)}^{1}(\mathbb{Z}/4\mathbb{Z}). \\
\end{eqnarray*}

\noindent v) Let $n=2l \geq 4$, $l$ even and $m \geq 3$. By lemma \ref{41} and lemma \ref{42}, we have:
\begin{eqnarray*}
w_{n+2}^{\overline{-1}}(\mathbb{Z}/2^{m}\mathbb{Z}) &=& 2^{n-1}w_{n+2}^{-1}(\mathbb{Z}/2^{m-1}\mathbb{Z})+2^{m(l-1)}\left|R_{2}^{\overline{-1}}(\mathbb{Z}/2^{m}\mathbb{Z})\right|-2^{m(l-1)-1}\sum_{k=0}^{p-1} \left|R_{2}^{\overline{-1+(p-k)p^{m-1}}}(\mathbb{Z}/2^{m} \mathbb{Z})\right| \\
                                                    &=& 2^{n-1}w_{n+2}^{-1}(\mathbb{Z}/2^{m-1}\mathbb{Z})+2^{m(l-1)}2^{m-1}(2(m+1)-m)-2^{m(l-1)-1}2^{m}(m(2-1)+1) \\
																									  &=& 2^{n-1}w_{n+2}^{-1}(\mathbb{Z}/2^{m-1}\mathbb{Z})+2^{ml-1}((2(m+1)-m)-(m(2-1)+1)) \\
																									  &=& 2^{n-1}w_{n+2}^{-1}(\mathbb{Z}/2^{m-1}\mathbb{Z})+2^{ml-1} \\
																									  &=& \ldots~= 2^{(m-2)(n-1)}w_{2(l+1)}^{-1}(\mathbb{Z}/4\mathbb{Z})+\sum_{j=0}^{m-3} 2^{j(n-1-l)}2^{ml-1} \\
																									  &=& 2^{(m-2)(n-1)}w_{2(l+1)}^{-1}(\mathbb{Z}/4\mathbb{Z})+\frac{2^{(m-2)(n-1-l)}-1}{2^{n-1-l}-1}2^{ml-1}.
\end{eqnarray*}

\noindent If we combine these formulae with theorems \ref{26}, \ref{27} and \ref{28}, we obtain the formulae given in the theorem. Similar calculations give us the remaining formulae. Note that the formulae for $w_{4}^{\pm \overline{1}}(\mathbb{Z}/p^{m}\mathbb{Z})$ are given in lemma \ref{42}.

\end{proof}

\subsection{Numerical applications}

\noindent The table below gives some values of $w_{n}^{\overline{1}}(\mathbb{Z}/N\mathbb{Z})$.

\begin{center}
\begin{tabular}{|c|c|c|c|c|c|}
\hline
  \multicolumn{1}{|c|}{\backslashbox{$n$}{\vrule width 0pt height 1.25em$N$}} & $\mathbb{Z}/8\mathbb{Z}$ & $\mathbb{Z}/9\mathbb{Z}$ & $\mathbb{Z}/16\mathbb{Z}$ & $\mathbb{Z}/18\mathbb{Z}$ & $\mathbb{Z}/27\mathbb{Z}$   \rule[-7pt]{0pt}{18pt} \\
  \hline
  4   & 20 & 21 & 48 & 63 & 81    \rule[-7pt]{0pt}{18pt} \\
	\hline
	  5   & 80 & 90 & 320 & 450 &  810    \rule[-7pt]{0pt}{18pt} \\
	\hline
  6  & 640 & 702 & 5120 & 7722 &  18~954  \rule[-7pt]{0pt}{18pt} \\
	\hline
	 7   & 5376 & 7371 & 86~016 & 154~791 &  597~051    \rule[-7pt]{0pt}{18pt} \\
	\hline
  8  & 45~312 & 70~227 & 1~452~032 & 3~019~761 & 17~078~283    \rule[-7pt]{0pt}{18pt} \\
\hline
	
\end{tabular}
\end{center}

\vspace{0.5cm}
\noindent The table below gives some values of $w_{n}^{\overline{-1}}(\mathbb{Z}/N\mathbb{Z})$.
\vspace{0.5cm}

\begin{center}
\begin{tabular}{|c|c|c|c|c|c|c|}
\hline
  \multicolumn{1}{|c|}{\backslashbox{$n$}{\vrule width 0pt height 1.25em$N$}} & $\mathbb{Z}/8\mathbb{Z}$ & $\mathbb{Z}/9\mathbb{Z}$ & $\mathbb{Z}/16\mathbb{Z}$ & $\mathbb{Z}/18\mathbb{Z}$ & $\mathbb{Z}/27\mathbb{Z}$   \rule[-7pt]{0pt}{18pt} \\
  \hline
  4   & 8 & 6 & 16 & 18 & 18   \rule[-7pt]{0pt}{18pt} \\
	\hline
	  5   & 80 & 90 & 320 & 450 & 810    \rule[-7pt]{0pt}{18pt} \\
	\hline
  6  & 800 & 999 & 6528 & 10~989 & 27~459     \rule[-7pt]{0pt}{18pt} \\
	\hline
	  7   & 5376 & 7371 & 86~016 & 154~791 & 597~051    \rule[-7pt]{0pt}{18pt} \\
	\hline
  8  & 43~008 & 63~180 & 1~376~256 & 2~716~740 & 15~352~740    \rule[-7pt]{0pt}{18pt} \\
\hline
	
\end{tabular}
\end{center}

\vspace{0.5cm}
\noindent \textbf{Acknowledgements.} The author thanks an anonymous referee for the very helpful comments and the many valuable suggestions regarding the presentation and the language of this paper.


\begin{thebibliography}{1}

\bibitem{BP}
C. Bardavid, É. Pité,
{\it Structure des anneaux commutatifs finis},
Revue de la Filière Mathématiques, Vol. 130 no. 4, (2020).

\bibitem{BC}
B. Böhmler, M. Cuntz,
{\it Frieze patterns over finite commutative local rings},
Séminaire Lotharingien de Combinatoire, Vol. 91, (2025), 18 pp.

\bibitem{CO}
C. Conley, V. Ovsienko,
{\it Rotundus : triangulations, Chebyshev polynomials, and Pfaffians}, 
The Mathematical Intelligencer, Vol. 40 no. 3, (2018), pp 45-50.

\bibitem{CO2}
C. Conley, V. Ovsienko,
{\it Quiddities of polygon dissections and the Conway-Coxeter frieze equation}, 
Annali della Scuola Normale Superiore di Pisa, Vol. 24 no. 4, (2023), pp 2125-2170.

\bibitem{CH} 
M. Cuntz, T. Holm,
{\it Frieze patterns over integers and other subsets of the complex numbers}, 
J. Comb. Algebra., Vol. 3 no. 2, (2019), pp 153-188.

\bibitem{CM} 
M. Cuntz, F. Mabilat,
\textit{Comptage des quiddités sur les corps finis et sur quelques anneaux $\mathbb{Z}/N\mathbb{Z}$}, 
Annales de la Faculté des Sciences de Toulouse, Vol 34 no. 1, (2025), pp 75-105.

\bibitem{E} 
L. Euler.
{\it Introductio in analysin infinitorum, Tomus primus}. 
Lausannae : apud Marcum Michaelum Bousquet et socios, 1748.

\bibitem{F}
B. Fine, 
\textit{Classification of finite rings of order $p^{2}$},
Mathematics Magazine, Vol. 66 no. 4, (1993), pp 248-252.

\bibitem{M4}
F. Mabilat,
\textit{Solutions monomiales minimales irréductibles dans $SL_{2}(\mathbb{Z}/p^{n}\mathbb{Z})$}, 
Bulletin des Sciences Mathématiques, Vol. 194, Article 103456, (2024), https://doi.org/10.1016/j.bulsci.2024.103456, hal-03573421, arxiv:2202.07279.

\bibitem{M5}
F. Mabilat,
\textit{Some counting formulae for $\lambda$-quiddities over the rings $\mathbb{Z}/2^{m}\mathbb{Z}$}, 
Bulletin of the Australian Mathematical Society, Vol. 111 no. 1, (2025), pp 1-12.

\bibitem{M6}
F. Mabilat,
\textit{Étude des liens entre la taille et l'irréductibilité des solutions monomiales minimales dans $SL_{2}(\mathbb{Z}/N\mathbb{Z})$}, 
(2024), arXiv:2406.18188.

\bibitem{Mo2}
S. Morier-Genoud,
\textit{Counting Coxeter's friezes over a finite field via moduli spaces}, 
Algebraic combinatoric, Vol. 4 no. 2, (2021), pp 225-240.

\end{thebibliography}
\end{document}